\documentclass[aos,MSNbibl,dvips]{arximspdf}
\usepackage{dcolumn}

\doi{10.1214/14-AOS1221} 
\volume{42}
\issue{3}
\pubyear{2014}
\firstpage{1166}
\lastpage{1202}

\makeatletter
\newcolumntype{d}[1]{D{.}{.}{#1}}
\newcommand{\rrvert}{\vert}
\newcommand{\llvert}{\vert}
\def\cal{\mathcal}
\renewcommand{\citep}[1]{(\citeauthor{#1} \citeyear{#1})}
\newcommand{\eqref}[1]{(\ref{#1})}
\newtheorem{theo}{Theorem}[section]

\newtheorem{lemm}{Lemma}[section]
\newtheorem{corr}{Corollary}[section]
\newproclaim{remark}{Remark}[section]
\newproclaim{example}{Example}[section]
\newproclaim{definition}{Definition}[section]
\newproclaim{condition}{Condition}[section]

\newcommand{\PP}{\mathbb{P}}
\newcommand{\EE}{\mathbb{E}}
\newcommand{\R}{\mathbb{R}}
\newcommand{\Nat}{\mathbb{N}}
\newcommand{\eps}{\varepsilon}
\newcommand{\Var}{\operatorname{Var}}
\newcommand{\bx}{\mathbf{X}}
\newcommand{\argmin}{\mathop{\arg\min}}

\def\summ#1#2#3{\sum_{#1=#2}^{#3}}

\makeatother

\begin{document}
\begin{frontmatter}

\title{On asymptotically optimal confidence regions and tests for
high-dimensional models}
\runtitle{Confidence regions for high-dimensional models}

\begin{aug}
\author[a]{\fnms{Sara}~\snm{van de Geer}\corref{}\ead[label=sg]{geer@stat.math.ethz.ch}},
\author[a]{\fnms{Peter}~\snm{B\"uhlmann}\ead[label=pb]{buhlmann@stat.math.ethz.ch}},
\author[b]{\fnms{Ya'acov}~\snm{Ritov}\ead[label=yr]{yaacov.ritov@gmail.com}\thanksref{t1}}
\and
\author[a]{\fnms{Ruben}~\snm{Dezeure}\ead[label=rd]{dezeure@stat.math.ethz.ch}}
\runauthor{van de Geer, B\"uhlmann, Ritov and Dezeure}
\affiliation{ETH Z\"urich, ETH Z\"urich, The Hebrew University of
Jerusalem and\break
ETH Z\"urich}
\thankstext{t1}{Supported by the
Forschungsinstitut f\"ur Mathe\-matik (FIM) at
ETH Z\"urich and from the Israel Science Foundation (ISF).}
\address[a]{S. van de Geer\\
P. B\"uhlmann\\
R. Dezeure\\
Seminar f\"ur Statistik\\
ETH Z\"urich \\
R\"{a}mistrasse 101\\
8092 Z\"{u}rich\\
Switzerland\\
\printead{sg}\\
\phantom{E-mail:\ }\printead*{pb}\\
\phantom{E-mail:\ }\printead*{rd}}
\address[b]{Y. Ritov\\
Department of Statistics\\
The Hebrew University of Jerusalem\\
Mount Scopus, Jeruzalem 91905\\
Israel\\
\printead{yr}} 
\end{aug}

\received{\smonth{3} \syear{2013}}
\revised{\smonth{1} \syear{2014}}

%
\begin{abstract}
We propose a general method for constructing confidence intervals and
statistical tests for single or low-dimensional components of a large
parameter vector in a high-dimensional model. It can be easily adjusted
for multiplicity taking dependence among tests into account. For linear models,
our method is essentially the same as in Zhang and Zhang
[\textit{J. R. Stat. Soc. Ser. B Stat. Methodol.} \textbf{76}
(2014) 217--242]: we analyze
its asymptotic properties and establish its asymptotic optimality in terms
of semiparametric efficiency. Our method
naturally extends to generalized linear models with convex loss
functions. We
develop the corresponding theory which includes a careful analysis for
Gaussian, sub-Gaussian and bounded correlated designs.
\end{abstract}

%

\end{frontmatter}

\section{Introduction}\label{sec1}

Much progress has been made over the last decade in high-dimensional
statistics where the number of unknown parameters greatly
exceeds sample size. The vast majority of work has been pursued for point
estimation such as consistency for prediction \cite{greenrit03,pb06}, oracle
inequalities and estimation of a high-dimensional parameter
\cite{cantao07,buneaetal06,zhang2008sparsity,geer07,MY08,brt09,negetal12,judetal11}
or variable selection \cite{mebu06,zhaoyu06,fanlv07,Wai08}. Other
references and exposition to a broad class of models can be found in
\cite{fanlv10} or \cite{pbvdg11}.

Very little work has been done for constructing confidence intervals,
statistical testing and assigning uncertainty in high-dimensional sparse
models. A major difficulty of the problem is the fact that sparse estimators
such as the lasso do not have a tractable limiting distribution:
already in
the low-dimensional setting, it depends on the unknown parameter
\cite{knfu00} and
the convergence to the limit is not uniform. Furthermore, bootstrap and
even subsampling techniques are plagued by noncontinuity of
limiting distributions. Nevertheless, in the low-dimensional setting, a modified
bootstrap scheme has been proposed; \cite{chatter11}
and \cite{chatter13} have recently proposed a residual based
bootstrap scheme. They provide consistency guarantees for the
high-dimensional setting; we consider this method in an empirical analysis
in Section~\ref{empiricalsection}.

Some approaches for quantifying uncertainty include the following. The
work in
\cite{WR08} implicitly contains the idea of sample splitting and corresponding
construction of $p$-values and confidence intervals, and the procedure
has been
improved by using multiple sample splitting and aggregation of dependent
$p$-values from multiple sample splits
\cite{memepb09}. Stability selection \cite{mebu10} and its modification
\cite{shah13} provides another route to estimate error measures for false
positive selections in general high-dimensional settings.
An alternative method for obtaining confidence sets is in the recent
work \cite{Meinshausen2013}.
From another and
mainly theoretical perspective, the work in \cite{judetal11}
presents necessary and sufficient conditions for recovery with the lasso
$\hat{\beta}$ in terms of $\|\hat{\beta} - \beta^0\|_{\infty}$, where
$\beta^0$ denotes the true parameter: bounds on
the latter, which hold with
probability at least say $1- \alpha$, could be used in
principle to construct (very) conservative confidence regions.
At a theoretical level, the paper~\cite{Nickl2013}
derives confidence intervals in $\ell_2$ for the case of two possible
sparsity levels. Other recent
work is discussed in Section~\ref{subsecrelwork} below.

We propose here a method which enjoys optimality properties when
making assumptions on the sparsity and design matrix of the model. For a
linear model, the procedure is as the one in
\cite{zhangzhang11}
and closely related to the method in
\cite{jamo13}. It is based on the lasso and is ``inverting'' the
corresponding KKT conditions. This yields a nonsparse estimator which has
a Gaussian (limiting) distribution. We show, within a sparse linear model
setting, that the estimator is optimal in the sense that it reaches the
semiparametric efficiency bound. The procedure can be used and is
analyzed for
high-dimensional sparse linear and generalized linear models and for
regression problems with general convex (robust) loss functions.

\subsection{Related work}\label{subsecrelwork}

Our work is closest to \cite{zhangzhang11} who proposed the
semiparametric approach for distributional inference in a high-dimensional
linear model. We take here a slightly different view-point, namely by
inverting the KKT conditions from the lasso, while relaxed projections are
used in \cite{zhangzhang11}. Furthermore, our paper extends the
results in
\cite{zhangzhang11} by: (i) treating generalized linear models and general
convex loss functions; (ii) for linear models, we give conditions under
which the procedure achieves the semiparametric efficiency bound and our
analysis allows for rather general Gaussian, sub-Gaussian and bounded
design. A related approach as in \cite{zhangzhang11} was proposed in
\cite{pb12} based on ridge regression which is clearly suboptimal and
inefficient with a detection rate (statistical power) larger
than $1/ \sqrt n$.

Recently, and developed independently, the work in \cite{jamo13}
provides a
detailed analysis for linear models
by considering a very similar procedure
as in \cite{zhangzhang11} and in our paper. They show that the detection
limit is indeed in the $1/\sqrt{n}$-range and they provide a minimax test
result; furthermore, they present extensive simulation results indicating
that the ridge-based method in
\cite{pb12} is overly conservative, which is in line with the theoretical
results. Their optimality results are interesting and are complementary to
the semiparametric optimality established here. Our results cover a
substantially broader range of non-Gaussian designs in
linear models, and we provide a
rigorous analysis for correlated designs with covariance matrix $\Sigma
\neq I$: the
SDL-test in \cite{jamo13} assumes that $\Sigma$ is known while we
carefully deal with the issue when $\Sigma^{-1}$ has to be estimated (and
arguing why, e.g., GLasso introduced in \cite{friedman2008sparse} is
not good for our purpose).
Another way and method to achieve distributional inference for
high-dimensional models is given in \cite{beletal11} (claiming
semiparametric efficiency). They use a two-stage
procedure with a so-called post-double-selection as first and least squares
estimation as second stage: as such, their methodology is radically
different from ours. At the time of writing of this paper,
\cite{JavanmardMontanari13} developed another modification
which directly computes an approximate inverse of the Gram matrix. Moreover,
\cite{Belloni2013a} extended their approach to logistic regression and
\cite{Belloni2013b} to LAD estimation using an instrumental variable approach.

\subsection{Organization of the paper}\label{subsecorganization}

In Section~\ref{seclinmod}, we consider the linear model and the
lasso. We describe the desparsifying step in Section~\ref
{subsecdespars} where we need to use an approximately
inverting matrix. A way to obtain this matrix is by applying the lasso
with nodewise regression, as given
in Section~\ref{subsecnodelasso}. Assuming Gaussian errors, we
represent in Section~\ref{subsecfixeddesign} the de-sparsified
lasso as sum of a normally distributed term and a remainder term.
Section~\ref{subsecranddesign} considers the case of random design
with i.i.d. covariables.
We first prove for the case of Gaussian design and Gaussian errors that
the remainder
term is negligible. We then show in Section~\ref{uniformsection} that
the results lead to honest
asymptotic confidence intervals. Section~\ref{subsecass} discusses the
assumptions
and Section~\ref{secasyEff} asymptotic efficiency.
The case of non-Gaussian design and
non-Gaussian errors is treated in Section~\ref{subsecsub-Gaussian}.

In Section~\ref{secGLM}, we consider the extension to generalized
linear models.
We start out in Section~\ref{desparseGLMsection}
with the procedure, which is again desparsifying the
$\ell_1$-penalized estimator. We again use the lasso with nodewise
regression to obtain an approximate inverse of the matrix
of second order derivatives. The computation of this approximate
inverse is
briefly described in Section~\ref{subsecnodelassomatrix}.
Section~\ref{subsectheoryGLM} presents asymptotic normality under
high-level conditions.
In Section~\ref{subsecnodewiserandommat}, we investigate the
consistency of
the lasso with nodewise regression as estimator of the inverse of the
matrix of second-order
derivatives of the theoretical risk evaluated at the true unknown
parameter $\beta^0$.
We also examine here the consistent estimation of the asymptotic variance.
Section~\ref{subsecassglm} gathers the results, leading to Theorem
\ref
{GLMtheorem} for generalized
linear models. Section~\ref{empiricalsection} presents some empirical results.
The proofs and theoretical material needed are given in Section~\ref
{secproofs}, while the technical proofs of
Section~\ref{secasyEff} (asymptotic efficiency)
and Section~\ref{subsecnodewiserandommat} (nodewise regression for
certain random matrices) are presented in the supplemental article
\cite{supplement}.

\section{High-dimensional linear models}\label{seclinmod}

Consider a high-dimensional linear model
%
\begin{equation}
\label{modlin} Y = \bx\beta^0 + \eps,
\end{equation}
with $n \times p$ design matrix $\bx=: [X_1,\ldots,X_p]$ ($n \times
1$ vectors $X_j$), $\eps\sim
{\cal N}_n(0,\sigma_{\eps}^2 I)$ independent of $\bx$ and unknown
regression $p \times1$ vector
$\beta^0$. We note that non-Gaussian errors are not a principal
difficulty, as
discussed in Section~\ref{subsecsub-Gaussian}. Throughout the paper,
we assume that $p > n$ and
in the asymptotic results we require $\log(p )/ n = o(1)$. We denote
by $S_0:= \{j; \beta_j^0 \neq0\}$ the active set of variables and its cardinality
by $s_0:= |S_0|$.

Our main goal is a pointwise statistical inference for the components of
the parameter vector $\beta^0_j\ (j=1,\ldots,p)$ but we also discuss
simultaneous inference for parameters $\beta^0_G:= \{\beta^0_j; j
\in
G\}$
where $G \subseteq\{1,\ldots,p\}$ is any group. To exemplify, we
might want to
test statistical
hypotheses of the form $H_{0,j}\dvtx \beta^0_j = 0$ or $H_{0,G}\dvtx \beta
^0_j =
0 \mbox{ for all } j \in G$, and when pursuing many tests, we aim for
an efficient multiple testing adjustment taking dependence into account and
being less conservative than say the Bonferroni--Holm procedure.

\subsection{The method: Desparsifying the lasso}\label{subsecdespars}

The main idea is to invert the Karush--Kuhn--Tucker characterization of the
lasso.

The lasso \cite{tibs96} is defined as
%
\begin{equation}
\label{lasso} \hat{\beta} = \hat{\beta}(\lambda):= \argmin_{\beta\in\R^p} \bigl(
\|Y - \bx \beta\|_2^2/n + 2 \lambda\|\beta
\|_1\bigr).
\end{equation}
It is well known that the estimator in (\ref{lasso}) fulfills the
Karush--Kuhn--Tucker (KKT) conditions:
\begin{eqnarray*}
& - \bx^T (Y - \bx\hat{\beta})/n + \lambda\hat{\kappa} = 0,&
\\
& \|\hat{\kappa}\|_{\infty} \le1\quad \mbox{and}\quad \hat{\kappa}_j
= \operatorname {sign}(\hat{\beta}_j)\qquad \mbox{if $\hat{\beta}_j
\neq0$}.&
\end{eqnarray*}
The vector $\hat{\kappa}$ is arising from the subdifferential of
$\|\beta\|_1$: using the first equation we can always represent it as
%
\begin{equation}
\label{hattau} \lambda\hat{\kappa} = \bx^T(Y - \bx\hat{\beta})/n.
\end{equation}
The KKT conditions can be rewritten with the notation $\hat{\Sigma} =
\bx^T \bx/n $:
\[
\hat{\Sigma} \bigl(\hat{\beta}- \beta^0\bigr) + \lambda\hat{\kappa}
= \bx^T \eps/n.
\]

The idea is now to use a ``relaxed form'' of an inverse of
$\hat{\Sigma}$. Suppose that $\hat{\Theta}$ is a reasonable approximation
for such an inverse, then
%
\begin{equation}
\label{pivot1} \hat{\beta} - \beta^0 +\hat{\Theta} \lambda\hat{\kappa}
= \hat{\Theta} \bx^T \eps/n - \Delta/\sqrt n,
\end{equation}
where
\[
\Delta:= \sqrt n (\hat{\Theta} \hat{\Sigma} - I) \bigl(\hat{\beta} -
\beta^0\bigr). %
\]
We will show in Theorem \ref{th2} that $\Delta$ is asymptotically
negligible under certain sparsity assumptions. This suggests the
following estimator:
%
\begin{equation}
\label{b1} \hat{b} = \hat{\beta} + \hat{\Theta} \lambda\hat{\kappa} = \hat{
\beta} + \hat{\Theta} \bx^T(Y - \bx\hat{\beta})/n,
\end{equation}
using (\ref{hattau}) in the second equation. This is essentially the same
estimator as in \cite{zhangzhang11} and it is of the same form as the
SDL-procedure in
\cite{jamo13}, when plugging in the estimate $\hat{\Theta}$ for the
population quantity $\Theta:= \Sigma^{-1}$ where $\Sigma$ is the
population inner
product matrix.
With (\ref{pivot1}), we immediately obtain an asymptotic pivot when
$\Delta$ is negligible, as is justified in Theorem
\ref{th2} below:
%
\begin{equation}
\label{piv1} \sqrt{n}\bigl(\hat{b} - \beta^0\bigr) =
W+o_{\PP}(1),\qquad  W|\bx\sim{\cal N}_p\bigl(0,
\sigma_{\eps}^2 \hat{\Theta} \hat{\Sigma} \hat{
\Theta}^T\bigr).
\end{equation}

An asymptotic pointwise confidence interval for $\beta^0_j$ is then
given by
\begin{eqnarray*}
& \bigl[\hat{b}_j - c(\alpha,n,\sigma_{\eps}),
\hat{b}_j + c(\alpha,n,\sigma_{\eps})\bigr],&
\\
& c(\alpha,n,\sigma_{\eps}): = \Phi^{-1}(1 - \alpha/2)
\sigma_{\eps} \sqrt{\bigl(\hat{\Theta} \hat{\Sigma} \hat{
\Theta}^T\bigr)_{j,j}/n},&
\end{eqnarray*}
where $\Phi(\cdot)$ denotes the c.d.f. of ${\cal N}(0,1)$.
If $\sigma_{\eps}$ is unknown, we replace it by a consistent estimator.

\subsubsection{The lasso for nodewise regression}\label{subsecnodelasso}

A prime example to construct the approximate inverse $\hat{\Theta}$ is
given by the lasso for the nodewise regression on
the design $\bx$: we use the lasso $p$ times for each
regression problem
$X_j$ versus $\bx_{-j}$, where the latter is the design submatrix without
the $j$th column. This method was introduced by \cite{mebu06}. We
provide here a
formulation suitable for our purposes.
For each $j = 1,\ldots,p$,
%
\begin{equation}
\label{nodelasso} \hat{\gamma}_j:= \argmin_{\gamma\in\R^{p-1} } \bigl(
\|X_j - \bx_{-j} \gamma\|_2^2/n + 2
\lambda_j \|\gamma\|_1\bigr),
\end{equation}
with components of $\hat{\gamma}_j =
\{\hat{\gamma}_{j,k}; k=1,\ldots,p, k \neq j\}$.
Denote by
\begin{eqnarray*}
\hat{C}:= \pmatrix{ 1 & -\hat\gamma_{1,2} & \cdots& - \hat
\gamma_{1,p} \vspace *{2pt}
\cr
-\hat \gamma_{2,1} & 1 & \cdots&
-\hat\gamma_{2,p}
\vspace*{2pt}\cr
\vdots& \vdots& \ddots& \vdots\vspace*{2pt}
\cr
- \hat\gamma_{p,1} & - \hat
\gamma_{p,2} & \cdots& 1 }
\end{eqnarray*}
and write
\[
\hat T^2:= \operatorname{diag} \bigl( \hat\tau_1^2, \ldots, \hat\tau_p^2 \bigr),
\]
where for $j=1, \ldots, p$
\[
\hat{\tau}_j^2:= \| X_j -
\bx_{-j} \hat{\gamma}_j \|_2^2 / n
+ \lambda_j \| \hat\gamma_j \|_1.
\]
Then define
%
\begin{equation}
\label{thetalasso} \hat{\Theta}_{\mathrm{Lasso}}: = \hat{T}^{-2} \hat{C}.
\end{equation}
Note that although $\hat\Sigma$ is self-adjoint, its relaxed inverse
$\hat{\Theta}_{\mathrm{Lasso}}$ is not. In the sequel, we denote by
%
\begin{equation}
\label{blasso} \hat{b}_{\mathrm{Lasso}}= \mbox{the estimator in (\ref{b1})
with $\hat \Theta$ the nodewise lasso from (\ref{thetalasso})}.
\end{equation}
The estimator $\hat{b}_{\mathrm{Lasso}}$ corresponds to the proposal
in \cite{zhangzhang11}.

Let the $j$th row of $\hat{\Theta}$ be denoted by
$\hat{\Theta}_j$ (as a $1 \times p$ vector) and analogously for~$\hat{C}_j$. Then $\hat{\Theta}_{\mathrm{Lasso},j} =
\hat{C}_j/\hat{\tau}_j^2$.

The KKT conditions for the nodewise lasso \eqref{nodelasso} imply that
\[
\hat\tau_j^2 = (X_j - \bx_{-j}
\hat{\gamma}_j)^T X_j/n
\]
so that
\[
X_j^T \bx\hat{\Theta}_{\mathrm{Lasso},j}^T /n =
1.
\]
These KKT conditions also imply that
\[
\bigl\| \bx_{-j}^T \bx\hat{\Theta}_{\mathrm{Lasso},j}^T
\bigr\|_{\infty} /n \le \lambda_j / \hat \tau_j^2.
\]
Hence, for the choice $\hat{\Theta}_j =\hat{\Theta}_{\mathrm{Lasso},j}
$ we have
%
\begin{equation}
\label{add-inver} \bigl\| \hat{\Sigma} \hat{\Theta}_j^T
-e_ j \bigr\|_{\infty} \le\lambda_j /\hat{
\tau}_j^2,
\end{equation}
where $e_j$ is the $j$th unit column vector.
We call this the extended KKT conditions.

We note that using, for example, the GLasso estimator of \cite
{friedman2008sparse} for $\hat{\Theta}$ may not be
optimal because with this choice a bound for $\| \hat\Sigma\hat
\Theta
_j^T - e_j \|_{\infty}$
is not readily available and this means we cannot directly derive
desirable componentwise properties of the estimator $\hat{b}$ in
(\ref{b1}) as established in Section~\ref{subsecranddesign}. The same
can be said about a ridge type of estimator for $\hat\Theta$, a choice
analyzed in \cite{pb12}.
We note that in (\ref{add-inver}) the bound depends on $\hat\tau_j^2$
and is in this sense not under control.
In~\cite{JavanmardMontanari13}, a program is proposed which gives an
approximate inverse $\hat\Theta$
such that $\| \hat\Sigma\hat\Theta_j^T - e_j \|_{\infty}$ is bounded
by a prescribed constant.
We will show in Remark \ref{studentizeremark} that a bound of the form
(\ref{add-inver}) with $\lambda_j$ proportional (by a prescribed
constant) to
$\tilde\tau_j:= \| X_j - \bx_{-j} \hat\gamma_j \|_2 /\sqrt n $ gives
the appropriate normalization when considering a Studentized
version of the estimator $\hat b_\mathrm{Lasso}$.


\subsection{Theoretical result for fixed design}\label{subsecfixeddesign}

We provide here a first result for fixed design $\bx$.
A crucial identifiability assumption on the design is the so-called
compatibility condition \cite{vandeGeer07a}.
To describe this condition, we introduce the following notation. For a
$p \times1$ vector $\beta$
and a subset $S \subseteq\{1,\ldots,p\}$, define $\beta_S$ by
\[
\beta_{S,j}:= \beta_j 1 \{ j \in S\},\qquad j=1,
\ldots,p.
\]
Thus, $\beta_S$ has zeroes for the components outside the set $S$. The
compatibility condition for $\hat\Sigma$ requires a positive constant
$\phi_0 > 0$ such that
for all $\beta$ satisfying $\|\beta_{S_0^c}\|_1 \le3 \|\beta_{S_0}\|
_1$ (the constant $3$
is relatively arbitrary, it depends on the choice of the tuning
parameter $\lambda$)
\[
\|\beta_{S_0}\|_1^2 \le{s_0}
\beta^T \hat{\Sigma} \beta/ \phi_0^2.
\]
The value $\phi_0^2$ is called the compatibility constant.

We make the
following assumption:
\begin{longlist}[(A1)]
\item[(A1)] The compatibility condition holds for $\hat{\Sigma}$ with
compatibility constant \mbox{$\phi_0^2> 0$}. Furthermore, $\max_j \hat
{\Sigma
}_{j,j} \le M^2
$ for some $0 < M < \infty$.
\end{longlist}
The assumption (A1) is briefly discussed in Section~\ref{subsecass}.
We then obtain the following result
where we use the notation $\|A\|_{\infty}:= \max_{j,k} |A_{j,k}|$ for
the element-wise sup-norm
for a matrix $A$.
%
\begin{theo}\label{th1}
Consider the linear model in (\ref{modlin}) with Gaussian error $\eps
\sim
{\cal N}_n(0,\sigma_{\eps}^2 I)$, and assume \textup{(A1)}. Let $t >0$ be arbitrary.
When using the lasso in
(\ref{lasso}) with $\lambda\ge2 M \sigma_{\eps} \sqrt{2(t^2 +
\log(p))/{n}}$ and
the lasso for nodewise regression in (\ref{thetalasso}) we have:
\begin{eqnarray*}
& \sqrt n \bigl( \hat{b}_{\mathrm{Lasso}} - \beta^0 \bigr) = W +
\Delta,&
\\
&  W= \hat{\Theta}_{\mathrm{Lasso}} \bx^T \eps/\sqrt n \sim{\cal
N}_n \bigl(0, \sigma_{\varepsilon}^2 \hat\Omega\bigr),\qquad
\hat\Omega:= \hat\Theta\hat\Sigma\hat\Theta^T,&
\\
&\displaystyle \PP \biggl[\|\Delta\|_{\infty} \ge8 \sqrt n \biggl( \max
_j \frac
{\lambda_j }{ \hat\tau_j^2} \biggr) \frac{\lambda s_0}{\phi_0^2 } \biggr] \le2
\exp\bigl[-{t^2}\bigr].&
\end{eqnarray*}
\end{theo}
A proof is given in Section~\ref{subsecproofth1}.

\begin{remark} \label{studentizeremark} In practice, one will use
a Studentized version of $\hat b_\mathrm{Lasso} $. Let us consider the
$j$th component.
One may verify that
$\hat\Omega_{j,j} =\tilde\tau_j^2 / \hat\tau_j^4 $,
where $\tilde\tau_j^2$ is the residual sum of squares
$\tilde\tau_j^2:= \| X_j - \bx_{-j} \hat\gamma\|_2^2 / n $.
Under the conditions of Theorem~\ref{th1},
\begin{eqnarray*}
&\displaystyle \frac{\sqrt n ( \hat{b}_{\mathrm{Lasso},j } - \beta_j^0 ) } {
\hat\Omega_{j,j}^{1/2} \sigma_{\varepsilon} } = V_j + \tilde\Delta _j,&
\\
& V_j \sim{\cal N} (0, 1),&
\\
&\displaystyle \PP \biggl[|\tilde\Delta_j | \ge8 \sqrt n \biggl(
\frac{\lambda_j
} { \tilde\tau_j } \biggr) \biggl(\frac{\lambda} { \sigma_{\varepsilon}} \biggr) \frac{ s_0}{\phi_0^2 }
\biggr] \le2 \exp\bigl[-{t^2}\bigr].&
\end{eqnarray*}
A Studentized version has the unknown variance $\sigma_{\varepsilon}^2$
replaced by a consistent estimator, $\tilde\sigma_{\varepsilon}^2$ say.
Thus, the bound for $\tilde\Delta_j $ depends on the normalized
tuning parameters $\lambda_j/ \tilde\tau_j$ and $\lambda/ \tilde
\sigma
_{\varepsilon}$.
In other words, the standardized estimator is standard normal with a
standardized remainder term.
The appropriate choice for $\lambda$ makes $\lambda/ \tilde\sigma
_{\varepsilon}$
scale independent. Scale independence for $\lambda_j/ \tilde\tau_j$
can be shown under
certain conditions, as we will do in the next subsection. Scale
independent regularization
can also be achieved numerically by
using the square-root lasso introduced in \cite{belloni2011square},
giving an approximate inverse, $\hat\Theta_{\sqrt\mathrm{Lasso}}$ say,
as alternative for $\hat\Theta_\mathrm{Lasso}$. Most of the theory that
we develop
in the coming subsections goes through with the choice $\hat\Theta
_{\sqrt\mathrm{Lasso}}$ as well. To avoid digressions, we do
not elaborate on this.
\end{remark}

Theorem \ref{th2} presents conditions that ensure that $\hat\tau_j$ as
well as $1/ \hat\tau_j^2$ are asymptotically bounded uniformly in $j$
(see Lemma \ref{lem3} in Section~\ref{secproofs}) and that
asymptotically one may choose $\lambda$
as well as each
$\lambda_j $ of order $ \sqrt{\log(p)/n}$.
Then, if the sparsity $s_0$ satisfies
$s_0=o(\sqrt n / \log p )$,
the correct
normalization
factor for $\hat{b}_{\mathrm{Lasso}}$ is $\sqrt{n}$ (as used in the
above theorem) and the error term
$\|\Delta\|_{\infty} = o_{\PP}(1)$ is negligible. The
details are discussed next.

\subsection{Random design and
optimality}\label{subsecranddesign}

In order to further analyze the error term $\Delta$ from Theorem
\ref{th1}, we consider an asymptotic
framework with random design.
It uses a scheme where $p = p_n \ge n \to\infty$ in
model (\ref{modlin}), and thus, $Y = Y_n$, $\bx= \bx_n$, $\beta^0 =
\beta^0_n$ and $\sigma^2_{\eps} = \sigma^2_{\eps,n}$ are all (potentially)
depending on $n$. In the sequel, we usually suppress the index
$n$. We make the following assumption.
\begin{longlist}[(A2)]
\item[(A2)] The rows of $\bx$ are
i.i.d. realizations from a Gaussian distribution whose
$p$-dimensional inner product matrix $\Sigma$ has strictly positive
smallest eigenvalue
$\Lambda^2_{\mathrm{min}} $ satisfying $1/ \Lambda_\mathrm{min}^2 =
{\mathcal O} (1) $.
Furthermore, $\max_j \Sigma_{j,j} = {\mathcal O}(1)$.
\end{longlist}
The Gaussian assumption is relaxed in Section~\ref{subsecsub-Gaussian}.

We will assume below sparsity with respect to rows of $\Theta:=
\Sigma^{-1}$ and define
\[
s_{j}: = \bigl| \{ k \neq j \dvtx \Theta_{j, k} \neq0 \} \bigr|.
\]
%

Recall the notation
$\hat\Omega:=\hat{\Theta}_{\mathrm{Lasso}} \hat{\Sigma} \hat
{\Theta
}_{\mathrm{Lasso}}^T$.
We then have the following main result.
\begin{theo}\label{th2}
Consider the linear model (\ref{modlin}) with Gaussian error $\eps
\sim
{\cal N}_n(0,\sigma_{\eps}^2 I)$ where $\sigma_{\varepsilon}^2 =
{\mathcal
O} (1)$. Assume \textup{(A2)} and the sparsity
assumptions $s_{0} = o(\sqrt n / \log(p))$ and $\max_j s_{j}
= o(n/\log(p))$.\vadjust{\goodbreak} Consider a suitable choice of the
regularization parameters
$\lambda\asymp\sqrt{\log(p)/n}$ for the lasso in
(\ref{lasso}) and $\lambda_j \asymp
\sqrt{\log(p)/n}$ uniformly in $j$ for
the lasso for nodewise regression in (\ref{thetalasso}). Then
\begin{eqnarray*}
\sqrt{n}\bigl(\hat{b}_{\mathrm{Lasso}} - \beta^0\bigr)& =& W +
\Delta,
\\
W|\bx&\sim&{\cal N}_{p}\bigl(0,\sigma^2_{\eps}
\hat\Omega\bigr),
\\
\|\Delta\|_{\infty} &=&o_{\PP}(1).
\end{eqnarray*}
%
Furthermore, $\|\hat\Omega- \Sigma^{-1}\|_{\infty} = o_{\PP}(1)$.
\end{theo}
A proof is given in Section~\ref{subsecproofth2}.

Theorem \ref{th2} has various implications. For a one-dimensional component
$\beta_j^0$ (with $j$ fixed),
we obtain for all $z \in\R$
%
\begin{equation}
\label{one-dimens} \PP \biggl[{ \sqrt{n}(\hat{b}_{\mathrm{Lasso};j} - \beta^0_j)
\over
\sigma_{\eps}
\sqrt{\hat\Omega_{j,j}}} \le z \Big|\bx \biggr] - \Phi(z) =
o_{\PP}(1).
\end{equation}
Furthermore, for any fixed group $G \subseteq\{1,\ldots,p\}$ which is
potentially large, we have that
for all $z \in\R$
\[
\PP \biggl[\max_{j \in G} { \sqrt{n}|\hat{b}_{\mathrm{Lasso};j} -
\beta
^0_j| \over
\sigma_{\eps} \sqrt{\hat\Omega_{j,j}}} \le z \Big| \bx
\biggr] - \PP \biggl[\max_{j \in G} { |W_{j}| \over
\sigma_{\eps} \sqrt{\hat\Omega_{j,j}}}\le z \Big|
\bx \biggr] = o_{\PP}(1).
\]
Therefore, conditionally on $\bx$, the asymptotic distribution of
\[
\max_{j \in G} {n}|\hat{b}_{\mathrm{Lasso};j}|^2 /
\sigma_{\eps}^2 {\hat\Omega_{j,j}} %
\]
under the null-hypothesis
$H_{0,G}; \beta^0_j = 0\ \forall j \in G$ is asymptotically
equal to the maximum of dependent $\chi^2 (1)$ variables $\max_{j \in G}
|W_{j}|^2/ \sigma_{\eps}^2 {\hat\Omega_{j,j}} $ whose distribution can
be easily simulated since $\hat\Omega$
is known. The unknown $\sigma_{\varepsilon}^2$ may be replaced by a
consistent estimator. For example, the scaled lasso \cite{sunzhang11} yields
a consistent estimator for $\sigma_{\varepsilon}^2$ under the assumptions
made for
Theorem \ref{th2}.


Theorem \ref{th2} is extended in Theorem \ref{nodetheorem} to the case
of non-Gaussian errors and non-Gaussian
design.

\subsubsection{Uniform convergence}\label{uniformsection}

The statements of Theorem \ref{th2} also hold in a uniform sense, and thus
the confidence intervals and tests based on these statements
are honest \cite{li89}. In particular, the estimator $\hat b_\mathrm{Lasso}$
does not suffer the problems arising from the nonuniformity of limit
theory for penalized estimators
(described in, e.g., \cite{potscher2009confidence} or \cite
{potscher2009adistribution}).
Such uniformity problems are also taken care of in~\cite
{berk2013valid} using an alternative procedure.
However, using $\hat b_\mathrm{Lasso}- \beta^0$ as pivot is
asymptotically less
conservative in general.

We consider the set of parameters
\[
{\cal B}(s) = \bigl\{\beta\in\R^p; \bigl| \{ j \dvtx \beta_j
\neq0 \} \bigr| \le s\bigr\}.
\]
We let $\PP_{\beta^0}$ be the distribution of the data under the linear
model (\ref{modlin}).
Then the following for $\hat{b}_{\mathrm{Lasso}}$ in (\ref{blasso}) holds.
%
\begin{corr}\label{cor1}
Consider the linear model (\ref{modlin}) with Gaussian error $\eps
\sim
{\cal N}_n(0,\sigma_{\eps}^2 I)$ where $\sigma_{\varepsilon}^2 =
{\mathcal
O} (1)$. Assume \textup{(A2)} and the sparsity
assumption $\beta^0 \in{\cal B} (s_0)$ with $s_{0} = o(\sqrt n / \log(p))$.
Suppose that $\max_j s_{j}
= o(n/\log(p))$. Then, when using suitable choices with $\lambda
\asymp
\sqrt{\log(p)/n}$ for the lasso in (\ref{lasso}), and $\lambda_j
\asymp\sqrt{\log(p)/n}$ uniformly $j$ for the lasso for nodewise
regression in (\ref{thetalasso})
\begin{eqnarray*}
& \sqrt{n}\bigl(\hat{b}_{\mathrm{Lasso}} - \beta^0\bigr) = W +
\Delta,&
\\
& W|\bx\sim{\cal N}_{p}\bigl(0,\sigma^2_{\eps}
\hat\Omega\bigr),\qquad \hat \Omega:= \hat\Theta\hat\Sigma\hat\Theta^T,&
\\
& \| \Delta\|_{\infty} =o_{\PP_{\beta^0} } (1) \qquad\mathrm{uniformly\ in\ }
\beta^0 \in{\cal B} (s_0).&
\end{eqnarray*}
%
Moreover, since $\hat\Omega$ does not depend
on $\beta^0$ we have as in Theorem
\ref{th2}, $\|\hat\Omega- \Sigma^{-1}\|_{\infty} = o_{\PP}(1)$.
\end{corr}

The proof is exactly the same as for Theorem \ref{th2} by simply noting
that $ \|\hat{\beta} - \beta^0\|_1 =
{\mathcal O}_{\PP_{\beta^0} }(s_{0} \sqrt{\log(p)/n})$ uniformly in
$\beta^0 \in{\cal B} (s_0)$ [with high probability, the compatibility
constant is bounded away from zero uniformly in all subsets $S_0$ with
$|S_0 | =o(\sqrt{n /\log(p)}$)].

Corollary \ref{cor1} implies that for $j \in\{1,\ldots,p\}$ and all
$z \in\R$,
\[
\sup_{\beta^0 \in{\cal B}(s_0 )}\biggl\llvert \PP_{\beta^0} \biggl[
{
\sqrt
{n}(\hat{b}_{\mathrm{Lasso};j} -
\beta_{j}^0) \over\sigma_{\eps} \sqrt{\hat\Omega_{j,j}}} \le z \Big|\bx \biggr] - \Phi(z) \biggr\rrvert =
o_{\PP} (1).
\]
%

Thus one can construct $p$-values for each component.
Based on many single $p$-values, we can use standard procedures for multiple
testing adjustment to control for various type I error measures. The
representation from Theorems \ref{th1} or~\ref{th2} with
$\|\Delta\|_{\infty}$ being sufficiently small allows to construct a
multiple testing adjustment which takes the dependence in terms of the
covariance $\hat\Omega$ (see Theorem~\ref{th2}) into account: the exact
procedure is described in \cite{pb12}. Especially when having strong
dependence among the $p$-values, the method is much less conservative than
the Bonferroni--Holm procedure for strongly controlling the family-wise error
rate.

\subsubsection{Discussion of the assumptions}\label{subsecass}

The compatibility condition in (A1) is weaker than many others which
have been proposed such as assumptions on restricted or sparse eigenvalues
\cite{van2009conditions}: a relaxation by a constant factor has
recently been given in \cite{sunzhang11}.
Assumption (A2) is rather weak in the sense that it
concerns the population inner product matrix. It implies condition (A1) with
$1/\phi_0={\mathcal O} (1)$ (see Lemma~\ref{lem2}) and $M = {\mathcal O} (1)$.

Regarding the sparsity assumption for $s_{0}$ in Theorem \ref{th1}, our
technique
crucially uses the $\ell_1$-norm bound $\|\hat{\beta} - \beta^0\|_1 =
{\mathcal O}_{\PP}(s_{0} \sqrt{\log(p)/n})$; see Lem\-ma~\ref{lem1}. In
order that this
$\ell_1$-norm converges to zero, the sparsity constraint
$s_{0} = o(\sqrt{n/\log(p)})$ is usually required. Our sparsity
assumption is slightly stricter by the factor $1/ \sqrt{\log(p)}$
(because the
normalization factor is $\sqrt{n}$), namely
$s_{0} = o(\sqrt n / \log(p))$.

\subsubsection{Optimality and semiparametric efficiency}\label{secasyEff}

$\!\!\!$Corollary \ref{cor1} establishes, in fact, that for any $j$, $\hat
b_{\mathrm{Lasso},j}$ is an asymptotically efficient estimator of
$\beta
^0_j$, in the sense that it is asymptotically normal with asymptotic
variance converging, as $n\to\infty$ to the variance of the best
estimator. Consider, the \emph{one}-dimensional sub-model,
%
\begin{equation}
\label{modsemip1} Y=\beta_j^0(X_j-
\bx_{-j}\gamma_j)+\bx_{-j}\bigl(
\beta^0_{-j}+\beta ^0_j\bx
_{-j}\gamma_j\bigr)+\eps,
\end{equation}
%
where
$X_j-\bx_{-j}\gamma_j$ is the projection in $L_2 (\PP)$ of $X_j$ to
the subspace
orthogonal to $\bx_{-j}$. Clearly, this is a linear submodel of the general
model \eqref{modlin}, passing through the true point. The Gauss--Markov
theorem argues that the best variance of an unbiased estimator of
$\beta_j^0$
in \eqref{modsemip1} is given by $\sigma_\eps^2
/(n \Var(X_{1,j}-\bx_{1, -j}\gamma_j))$. Corollary $\ref{cor1}$
shows that
$\sigma_\eps^2
/ \Var(X_{1,j}-\bx_{1, -j}\gamma_j)$
this is the asymptotic variance of $\sqrt n ( \hat b_{\mathrm
{Lasso},j}- \beta_j^0) $. Thus, $\sqrt n (\hat
b_{\mathrm{Lasso},j}- \beta_j^0)$ is asymptotically normal, with the
variance of the best
possible unbiased estimator. Note, that any regular estimator (regular at
least on parametric sub-models) must be asymptotically unbiased.

The main difference between this and most of the other papers on complex
models is that usually the lasso is considered as solving a nonparametric
model with parameter whose dimension $p$ is increasing to infinity, while
we consider the
problem as a semiparametric model in which we concentrate on a
low-dimensional model of interest, for example, $\beta^0_j$, while the
rest of the
parameters, $\beta^0_{-j}$, are considered as nuisance parameters.
That is,
we consider the problem as a semiparametric one.

In the rest of this
discussion, we put the model in a standard semiparametric framework
in which there is an infinite-dimensional population model.
Without loss of generality, the parameter of interest is $\beta^0_1$,
that is, the first component (extension to more than one but finitely many
parameters of interest is straightforward). Consider the random design
model where
the sequence $\{(Y_i, X_{i,1}, Z_i ) \}_{i=1}^{\infty}$ is i.i.d. with
%
\begin{equation}
\label{modsemip2} Y_1 = \beta^0_1X_{1,1}+K(Z_1)
+ \eps_1,\qquad \eps_1 \sim{\cal N}\bigl(0,
\sigma^2_{\eps}\bigr),
\end{equation}
where $\beta_1^0 \in\R$ is an unknown parameter and $K (\cdot)$ is an
unknown function.
When observing $\{(Y_i, X_{i,1}, Z_i ) \}_{i=1}^n$ this is the
partially linear regression model, where $\sqrt n$-consistency for
the parametric part $\beta_1^0$ can be achieved \cite{robinson1988root}.
We
observe the i.i.d. sequence $\{ (Y_i, X_{i,1}, \{ X_{i,j}^n \}
_{j=2}^{p_n} \} )\}_{i=1}^n $
such that
%
\begin{eqnarray}
\label{asyEff} \qquad&& Y_1 = \beta_1^0
X_{1,1}+ \summ j2{p_n} \beta^{n}_j
X_{1,j}^n+\eps _1^n,
 \nonumber\\
\eqntext{\eps_1^n\mbox{ independent of } X_{1,1},
X_{1,2}^n, \ldots,X_{1,p_n}^n,}
\\
&&
\EE \biggl[ K(Z_1)-\sum_{j \in S_n \cap\{ 2, \ldots,
p_n\} }
\beta^n_jX_{1,j}^n \biggr]
^2 \rightarrow0, \qquad|S_n | = o\bigl(\sqrt n / \log(p)\bigr),
\\
&& \EE \Biggl[ \EE[X_{1,1}|Z_1]-\summ j2{p_n}
\gamma^n_{1,j} X_{1,j}^n
\Biggr]^2 \rightarrow0,\nonumber
\\
&&\Biggl(K(Z_1)-\summ j2{p_n} \beta^n_jX_{1,j}^n
\Biggr) \Biggl(\EE[X_{1,1}|Z_1]-\summ j2{p_n}
\gamma^n_{1,j}X_{1,j}^n \Biggr)=
o_{\PP}\bigl(n^{-1/2}\bigr).\nonumber
\end{eqnarray}

\begin{theo}\label{thspj}
Suppose \eqref{asyEff} and the conditions of Theorem \ref{th2} are
satisfied, then
\[
\hat{b}_{\mathrm{Lasso};1}=\beta^0_1+\frac{1}n
\summ i1n \bigl(X_{i,1}-\EE[X_{i,1}|Z_i] \bigr)
\eps_i+o_{\PP}\bigl(n^{-1/2}\bigr).
\]
In particular, the limiting variance of
$\sqrt{n}(\hat{b}_{\mathrm{Lasso};1} - \beta^0_1)$ reaches the information
bound $\sigma^2_{\eps}/\EE(X_{1,1} - \EE[X_{1,1}|Z_1])^2$. Furthermore,
$\hat{b}_{\mathrm{Lasso};1}$ is regular at the one-dimensional parametric
sub-model with component $\beta_1^0$, and hence, $\hat{b}_{\mathrm{Lasso};1}$
is asymptotically efficient for estimating $\beta_1^0$.
\end{theo}
A proof is given in the supplemental article \cite{supplement}.

As a concrete example consider the following situation:
%
\begin{eqnarray}
\label{modsemip}   K(Z_1)& =& \sum_{j=2}^{\infty}
\beta^0_j X_{1,j}, 
\nonumber
\\[-8pt]
\\[-8pt]
\nonumber
 X_{1,j}^n&\equiv& X_{1,j} \qquad\forall j=1,
\ldots,p_n,
\end{eqnarray}
\mbox{where}
\begin{eqnarray}
&&\beta^0 \in{\cal B}
(s_0):= \bigl\{(\beta _j)_{j \in\Nat}; | \{ j \dvtx
\beta_j \neq0 \} | \le s_0\bigr\},\nonumber\\
\eqntext{ s_0 <\infty\mbox{ fixed},
\EE[X_{1,j} ] = 0\ \forall j \mbox{ and }\displaystyle \max_{j \in
\Nat
} \operatorname{var} ( X_{1,j} ) <
\infty,}
\\
&&\min_{S \subset\Nat} \Lambda_{\mathrm{min}}^2(S)>0,
\nonumber\\
\eqntext{ \mbox{where } \Lambda_{\mathrm{min}}^2(S) \mbox{ is
the smallest eigenvalue of the}}\\
\eqntext{\mbox{covariance matrix of } \{X_{1, j}\dvtx j \in S\},}
\\
& & \bigl| \{ k\dvtx \gamma_{1,k}\neq0 \}
\bigr| < \infty,
\nonumber\\
\eqntext{\mbox{where }\displaystyle \gamma_1:= \argmin_{\gamma\in\R_{\infty} }
\EE _{\PP
}\Biggl[\Biggl(X_{1,1} - \sum
_{k =2}^{\infty} \gamma_k X_{1, k}
\Biggr)^2\Biggr].}
\end{eqnarray}
%

Note that the assumption about the minimal eigenvalues $\{ \Lambda
_{\mathrm{min}}^2 (S)\dvtx S \subset\Nat\}$
is equivalent to saying that $\{ X_{1,j}\} _{j \in\Nat}$ has a positive
definite covariance function.

\begin{lemm}\label{concretelem}
Condition \eqref{asyEff} is satisfied in the above example.
\end{lemm}
A proof of this lemma is given in the supplemental article \cite{supplement}.

\subsubsection{Non-Gaussian design and non-Gaussian errors}\label
{subsecsub-Gaussian}

We extend Theorem~\ref{th2} to allow for non-Gaussian designs and
non-Gaussian errors.
Besides
covering
a broader range for linear models, the result is important for
the treatment of generalized linear models in Section~\ref{secGLM}.

Consider a random design matrix $\bx$ with i.i.d. rows having inner
product matrix $\Sigma$ with its inverse (assumed to exist) $\Theta=
\Sigma^{-1}$. For
$j=1, \ldots, p$, denote by
$\gamma_j:= \argmin_{\gamma\in\R^{p-1} } \EE[\|X_j - \bx_{-j}
\gamma
\|_2^2]$.
Define the
error $\eta_j:= X_j - \bx_{-j} \gamma_j$ with variance $\tau_j^2 =
\EE[\|\eta_j\|_2^2/n] = 1/\Theta_{j,j}$, $j=1, \ldots, p$. We make
the following assumptions:
\begin{longlist}[(B1)]
\item[(B1)] The design $\bx$ has either
i.i.d. sub-Gaussian rows (i.e.,\break
$\max_i \sup_{\| v \|_2 \le1 } \EE \exp[ |\sum_{j=1}^p v_j \mathbf
{X}_{i,j}| ^2 / L^2 ]
= {\mathcal O} (1) $ for some fixed constant\break  $L >0$)
or i.i.d. rows and for some $K \ge1$, $\|\bx
\|_{\infty} = \max_{i,j}|\bx_{i,j}| = {\mathcal O}(K)$.
The latter we call the bounded case. The strongly bounded case assumes in
addition that $\max_j \| \bx_{-j} \gamma_j \|_{\infty} = {\mathcal
O}(K) $.

\item[(B2)] In the
sub-Gaussian case, it holds that
$ \max_j \sqrt{s_j \log(p)/ n } = o(1)$. In the (strongly) bounded
case, we assume
that $\max_j K^2 s_j \sqrt{\log(p)/ n } = o (1)$.

\item[(B3)] The smallest eigenvalue $\Lambda_{\mathrm{min}}^2$ of
$\Sigma$ is strictly
positive and $1/ \Lambda_{\mathrm{min}}^2 = {\mathcal O} (1) $.
Moreover, $\max_{j} \Sigma_{j,j} = {\mathcal O}(1)$.
%
\item[(B4)] In the bounded case, it holds that $\max_j \EE\eta
_{1,j}^{4} = {\mathcal O} (K^4)$.
\end{longlist}

We note that the strongly bounded case in (B1) follows from the
bounded case if $\| \gamma_j \|_1 = {\mathcal O} (1) $.
Assumption (B2) is a standard sparsity
assumption for $\Theta$. Finally, assumption (B3) implies that $\|
\Theta_j \|_2 \le\Lambda_{\mathrm{min}}^{-2} = {\mathcal
O} (1)$ uniformly in $j$ so that in particular
$\tau_j^2 = 1/ \Theta_{j,j}$ stays away from zero.
Note that (B3) also implies $\tau_j^2 \le\Sigma_{j,j} = {\mathcal
O}(1) $ uniformly in $j$.

To streamline the statement of the results, we
write $K_0=1$ in the sub-Gaussian case and $K_0 = K$ in the (strongly)
bounded case.

\begin{theo} \label{nodetheorem} Suppose the conditions \textup{(B1)--(B4)}
hold. Denote by $\hat{\Theta}:= \hat\Theta_{\mathrm{Lasso}} $ and
$\hat{\tau}_j^2$, $j=1, \ldots, p$ the estimates from
the nodewise lasso in (\ref{thetalasso}). Then for suitable tuning
parameters $\lambda_j
\asymp K_0 \sqrt{ \log(p)/ n } $ uniformly in $j$, we have\looseness=-1
\begin{eqnarray*}
\| \hat\Theta_j - \Theta_j
\|_1& =& {\mathcal O}_{\PP} \biggl( K_0
s_j \sqrt{\log(p) \over n } \biggr), \\
\| \hat
\Theta_j - \Theta_j \|_2 &= &{\mathcal
O}_{\PP} \biggl( K_0 \sqrt{s_j \log(p)\over n } \biggr),
\\
\bigl| \hat\tau_j^2 -
\tau_j^2\bigr| &=& {\mathcal O}_{\PP} \biggl(
K_0 \sqrt{s_j \log(p) \over n } \biggr),\qquad j=1, \ldots, p.
\end{eqnarray*}\looseness=0
Furthermore,
%
\begin{eqnarray}
\bigl| \hat\Theta_j \Sigma\hat\Theta_j^T -
\Theta_{j,j} \bigr| \le \| \Sigma\|_{\infty} \| \hat
\Theta_j - \Theta_j \|_1^2
\wedge \Lambda _\mathrm{max}^2 \| \hat\Theta_j -
\Theta_j \|_2^2 + 2\bigl| \hat
\tau_j^2 - \tau_j^2\bigr|,\nonumber\\
 \eqntext{j=1,
\ldots, p,} %
\end{eqnarray}
where $\Lambda_\mathrm{max}^2$ is the maximal eigenvalue of $\Sigma$.
In the sub-Gaussian or
strongly bounded case the results are uniform in $j$.

Finally, assume model (\ref{modlin}) but assume instead of Gaussian
errors that $\{ \varepsilon_i \}_{i=1}^n$ are i.i.d. with
variance $\sigma_{\varepsilon}^2 = {\mathcal O} (1)$. Assume moreover in
the sub-Gaussian case for $\bx$ that the
errors are subexponential, that is, that\break $\EE\exp[ |\varepsilon_1 | /  L
] = {\mathcal O} (1)$
for some fixed $L$.
Apply the estimator (\ref{lasso}) with $\lambda\asymp K_0 \sqrt{
\log
(p) / n}$ suitably chosen. Assume that
$K_0 s_0 \log(p) / \sqrt n = o(1)$ and\break  $\max_j K_0 s_j \sqrt{\log
(p)/n} = o(1)$. Then
we have
\begin{eqnarray*}
& \sqrt{n}\bigl(\hat{b}_{\mathrm{Lasso}} - \beta^0\bigr) = W +
\Delta,&
\\
& W= \hat\Theta\bx^T \varepsilon/\sqrt n,&
\\
& | \Delta_j | =o_{\PP}(1) \qquad\forall j&
\end{eqnarray*}
and in the sub-Gaussian or strongly bounded case
\[
\| \Delta\|_{\infty} = o_{\PP}(1). %
\]
\end{theo}
A proof is given in Section~\ref{subsecprthenode}.

Note that the result is as in Theorem \ref{th2} except that $W \vert
\bx
$ is not necessarily normally
distributed. A
central limit theorem argument can be used to obtain approximate
Gaussianity of components of $W|\bx$ of fixed dimension.
This can also be done for moderately growing dimensions
(see, e.g., \cite{portnoy1987central}), which is useful for testing
with large groups $G$.

\section{Generalized linear models and general convex loss
functions}\label{secGLM}

We show here that the idea of de-sparsifying $\ell_1$-norm penalized
estimators and corresponding theory from Section~\ref{seclinmod}
carries over to models with convex loss functions such as generalized
linear models (GLMs).

\subsection{The setting and de-sparsifying the \texorpdfstring{$\ell_1$-norm}{l1-norm} regularized
estimator}\label{desparseGLMsection}

We consider the following framework with $1 \times p$ vectors of
covariables $x_i \in{\cal X}
\subseteq\R^p$ and univariate responses $y_i \in{\cal Y} \subseteq
\R$ for $i=1,\ldots,n$. As before, we denote by $\bx$ the design matrix
with $i$th row equal to $x_i$. At the moment, we do not distinguish
whether $\bx$ is random or fixed (e.g., when conditioning on $\bx$).

For $y \in{\cal Y}$ and
$x \in{\cal X}$ being a $1 \times p$ vector, we have a loss function
\[
\rho_{\beta} (y,x)= \rho(y, x \beta) \qquad\bigl(\beta\in\R^p
\bigr),
\]
which is assumed to be a strictly convex function in $\beta\in\R^p$.
We now define
\[
\dot\rho_{\beta}:= {\partial\over\partial\beta} \rho_{\beta}, \qquad\ddot
\rho_{\beta}:= {\partial\over\partial\beta\,\partial\beta^T } \rho_{\beta},
\]
where we implicitly assume that the derivatives exist.
For a function $g\dvtx {\cal Y} \times{\cal X} \to\R$,
we write\vspace*{1.5pt}
$P_n g:= \sum_{i=1}^n g(y_i,x_i) / n $ and $Pg:= \EE P_n g $.
Moreover, we let $\| g \|_n^2:= P_n g^2 $ and $\| g \|^2:= P g^2 $.

The $\ell_1$-norm regularized estimator is
%
\begin{equation}
\label{lasso-genrho} \hat{\beta} = \argmin_{\beta} \bigl(P_n
\rho_{\beta} + \lambda\| \beta\|_1\bigr).
\end{equation}
As in Section~\ref{subsecdespars}, we desparsify the estimator. For this
purpose, define
%
\begin{equation}
\label{hatSigmaglm} \hat\Sigma:= P_n \ddot\rho_{\hat\beta}.
\end{equation}
Note that in general, $\hat{\Sigma}$ depends on $\hat{\beta}$ (an exception
being the squared error loss). We construct $\hat\Theta=
\hat{\Theta}_{\mathrm{Lasso}}$ by doing a nodewise lasso with $\hat
\Sigma$
as input as detailed below in (\ref{nodelassoglm}).
We then define
%
\begin{equation}
\label{estglm} \hat{b}:= \hat\beta- \hat\Theta P_n \dot
\rho_{\hat\beta}.
\end{equation}
The estimator in (\ref{b1}) is a special case of (\ref{estglm}) with
squared error loss.

\subsubsection{Lasso for nodewise regression with matrix
input}\label{subsecnodelassomatrix}

Denote by $\hat{\Sigma}$ a matrix which we want to approximately invert
using the nodewise lasso. For every row $j$, we consider the optimization
%
\begin{equation}
\label{nodelasso2} \hat{\gamma}_j:= \argmin_{\gamma\in\R^{p-1} } \bigl(\hat{
\Sigma }_{j,j} - 2 \hat{\Sigma}_{j, \setminus j} \gamma+
\gamma^T \hat{\Sigma}_{\setminus j,\setminus j} \gamma+ 2 \lambda_j
\|\gamma\|_1\bigr),
\end{equation}
where $\hat{\Sigma}_{j, \setminus j}$ denotes the $j$th row of
$\hat{\Sigma}$ without the diagonal element $(j,j)$, and
$\hat{\Sigma}_{\setminus j,\setminus j}$ is the submatrix without the
$j$th row and $j$th column. We note that for the case where $\hat
{\Sigma}
= \bx^T \bx/n$, $\hat{\gamma}_j$ is the same as in (\ref{nodelasso}).

Based on $\hat{\gamma}_j$ from (\ref{nodelasso2}), we compute
%
\begin{equation}
\label{nodelasso2b} \hat{\tau}_j^2 = \hat{
\Sigma}_{j,j} - \hat{\Sigma}_{j,\setminus j} \hat{\gamma}_j.
\end{equation}
Having $\hat{\gamma}_j$ and $\hat{\tau}_j^2$ from (\ref
{nodelasso2}) and
(\ref{nodelasso2b}), we define the nodewise lasso as
%
\begin{eqnarray}
\label{nodelassoglm}  & \hat{\Theta}_{\mathrm{Lasso}} \mbox{ as in
(\ref{thetalasso}) using (\ref{nodelasso2})--(\ref{nodelasso2b}) from matrix input $
\hat{\Sigma}$ in (\ref{hatSigmaglm})}.
\end{eqnarray}
%
Moreover, we denote by
\[
\hat{b}_{\mathrm{Lasso}}:= \hat{b} \mbox{ from (\ref{estglm}) using the
nodewise lasso from (\ref{nodelassoglm})}.
\]

Computation of (\ref{nodelasso2}), and hence of $\hat{\Theta}$ can be
done efficiently via coordinate descent using the KKT conditions to
characterize the zeroes. Furthermore, an active set strategy leads to
additional speed-up. See, for example, \cite{friedetal09} and
\cite{meier06grouplasso}.

For standard GLMs, the matrix input $\hat{\Sigma} = P_n \ddot
\rho_{\hat\beta}$ in (\ref{hatSigmaglm}) can be written
as $\hat{\Sigma} = \hat\Sigma_{\hat{\beta}} = \bx_{{\hat{\beta}}}^T
\bx_{{\hat{\beta}}}/n$ with $\bx_{{\hat{\beta}}}:= W_{\beta}
\bx
$ and $W_{\hat{\beta}} = \operatorname{diag}(w_{\hat{\beta}})$ for some weights
$w_{i,\hat{\beta}} = w_{\hat{\beta}} (y_i, x_i)\ (i=1,\ldots,n)$.
Then we
can simply use the nodewise lasso as in (\ref{thetalasso}) but based
on the
design matrix
$\bx_{{\hat{\beta}}}$: in particular, we can use the standard lasso
algorithm.

\subsection{Theoretical results}\label{subsectheoryGLM}

We show here that the components of the estimator~$\hat{b}$ in
(\ref{estglm}), when normalized with the easily computable standard error,
converge to a standard Gaussian distribution. Based on such a result,
the construction of confidence intervals and tests is straightforward.

Let $\beta^0 \in\R^p$ be the unique minimizer of $P
\rho_{\beta}$ with $s_0$ denoting the number of nonzero coefficients.
We use analogous notation as in Section~\ref{subsecranddesign} but with
modifications for the current context.
The asymptotic framework, which allows for Gaussian approximation of
averages, is as in Section~\ref{subsecranddesign} for $p = p_n \ge n
\to
\infty$, and thus, $Y:= (y_1,\ldots,y_n )^T= Y_n$, $\bx= \bx_n$,
$\beta^0 =
\beta^0_n$ and underlying parameters are all (potentially) depending on
$n$. As before, we usually suppress the corresponding index $n$.


We make the following assumptions which are discussed in
Section~\ref{subsecassglm}. Thereby, we assume (C3), (C5), (C6) and (C8)
for some constant $K \ge1$ and positive constants $\lambda_*$ and $s_*$.
The constant $\lambda$ is the tuning parameter in (\ref{lasso-genrho}).
In Section~\ref{subsecassglm}, we will discuss the conditions
with $\lambda\asymp\sqrt{\log p / n }$ and for all $j$, $\lambda_*
\asymp\lambda_j \asymp
\sqrt{\log(p) / n } $ where $\lambda_j$ is the tuning parameter in
(\ref{nodelasso2}).
Moreover, there we will assume $s_* \ge s_j$ for all $j$. Here, $s_j =
| \{ k\neq j \dvtx \Theta_{\beta^0, j, k } \neq0 \} | $, $j=1, \ldots,
p$ with $\Theta_{\beta^0}:= ( P \ddot\rho_{\beta^0 } )^{-1}$ (assumed
to exist).
%
\begin{longlist}[(C1)]
\item[(C1)]
The derivatives
\[
\dot\rho(y, a ):= {d \over d a} \rho(y, a),\qquad \ddot\rho(y, a ):=
{ d^2 \over d a^2 } \rho(y,a), %
\]
exist for all $y, a$, and for some $\delta$-neighborhood ($\delta>0$),
$\ddot\rho( y, a )$ is
Lip\-schitz:
\[
\max_{a_0 \in\{ x_i \beta^0\}} \sup_{ |a- a_0 | \vee| \hat a - a_0|
\le\delta} \sup
_{y \in{\cal Y}} { | \ddot\rho(y, a) - \ddot\rho
(y, \hat a) | \over|a - \hat a | } \le1.
\]
Moreover,
\[
\max_{a_0 \in\{ x_i \beta^0\}} \sup_{y \in{\cal Y} } \bigl| \dot\rho(y,
a_0)\bigr| = {\mathcal O}(1), \qquad\max_{a_0 \in\{ x_i \beta^0\}} \sup
_{|a - a_0
| \le\delta} \sup_{y
\in{\cal Y} }\bigl | \ddot\rho(y, a)\bigr| = {
\mathcal O}(1). %
\]


\item[(C2)] It holds that
$\| \hat\beta- \beta^0 \|_1 = {\mathcal
O}_{\PP} ( s_0 \lambda)$, $\|\bx(\hat{\beta} - \beta^0)\|^2 =
{\mathcal O}_{\PP} ( s_0 \lambda^2 )$, and $\| \bx( \hat\beta-
\beta^0
) \|_n^2 = {\mathcal O}_{\PP} ( s_0 \lambda^2 )$.

\item[(C3)] It holds that $ \|\mathbf{X} \|_{\infty}:= \max_{i,j}
|\bx
_{i,j} | ={\mathcal O}( K) $.

\item[(C4)] It holds that
$ \| P_n \ddot\rho_{\hat\beta} \hat\Theta_j^T - e_j \|_{\infty} =
{\mathcal O}_{\PP}(\lambda_*) $.


\item[(C5)] It holds that
$\|\bx\hat\Theta_j^T \|_{\infty} = {\mathcal O}_{\PP} (K)$ and $\|
\hat
\Theta_j \|_1 = {\mathcal O}_{\PP} (\sqrt{s_*}) $.

\item[(C6)] It holds that
$ \| (P_n - P) \dot\rho_{\beta^0} \dot\rho_{\beta^0 }^T \|
_{\infty}=
{\mathcal O}_{\PP} (K^2 \lambda) $ and
moreover
\[
\max_j 1/\bigl(\hat\Theta P \dot\rho_{\beta^0} \dot
\rho_{\beta^0 }^T \hat \Theta^T\bigr)_{j,j}
= {\mathcal O} (1).
\]

\item[(C7)] For every $j$, the random variable
\[
{ \sqrt n ( \hat\Theta P_n \dot\rho_{\beta^0} )_j \over\sqrt{ (
\hat\Theta P \dot\rho_{\beta^0} \dot\rho_{\beta^0 }^T \hat
{\Theta}^T
)_{j,j} }} %
\]
converges weakly to a ${\cal N}(0,1)$-distribution.


\item[(C8)] It holds that
\[
K s_0 \lambda^2 =o \bigl(n^{-{1/2}} \bigr),\qquad
\lambda_* \lambda s_0 =o \bigl(n^{-{1/2}} \bigr) \quad\mbox{and}\quad
K^2 s_* \lambda+ K^2 \sqrt{s_0} \lambda= o(1).
\]
\end{longlist}
%
The following main result holds for fixed or random design according to
whether the assumptions hold for one or the other case.
%
\begin{theo} \label{th3}
Assume \textup{(C1)--(C8)}. For the estimator in (\ref{estglm}), we have for
each $j
\in\{1,\ldots,p\}$:
\[
\sqrt{n}\bigl(\hat{b}_j - \beta_{j}^0\bigr)/
\hat{\sigma}_j = V_j + o_\mathbf {P} (1),
\]
where $V_j $ converges weakly to a ${\cal N}(0,1)$-distribution and
where
\[
\hat\sigma_j^2:= \bigl( \hat\Theta P_n
\dot\rho_{\hat\beta} \dot \rho_{\hat\beta}^T \hat{
\Theta}^T \bigr)_{j,j}. %
\]
\end{theo}
A proof is given in Section~\ref{th3section}.
%
%
Assumption (C1) of Theorem \ref{th3} means that we regress to the classical
conditions for asymptotic normality in the one-dimensional case as in, for
example, \cite{cramer1946mathematical}. Assumption (C8) is a sparsity
assumption: for $K = O(1)$ and choosing $\lambda_* \asymp
\lambda\asymp\sqrt{\log(p)/n}$ the condition reads as $s_0 =
o(\sqrt{n}/\log(p))$ (as in Theorem \ref{th2}) and $s_* = o(\sqrt
{n/\log
(p)})$. All the other assumptions (C2)--(C7)
follow essentially from the conditions of Corollary \ref{node2corollary}
presented later, with the exception that (C3) is straightforward to
understand. For more details, see Section~\ref{subsecassglm}.


\subsection{About nodewise regression with certain random
matrices}\label{subsecnodewiserandommat}

We justify in this section most of the assumptions for Theorem
\ref{th3} when using the nodewise lasso estimator $\hat{\Theta} =
\hat{\Theta}_{\mathrm{Lasso}}$ as in (\ref{nodelassoglm})
and when\vadjust{\goodbreak} the matrix input is parameterized by $\hat{\beta}$ as for standard
generalized linear models. For notational simplicity, we drop the subscript
``lasso'' in $\hat{\Theta}$. Let
$w_{\beta} $ be an $n$-vector with entries
$w_{i, \beta} = w_{\beta} (y_i, x_i)$. We consider the
matrix $\bx_{{\beta}}:= W_{\beta} \bx$ where $W_{\beta} = \operatorname{diag}
(w_{\beta})$. We define $\hat\Sigma_{\beta}: =\bx_{\beta}^T
\bx_{{\beta} } / n $. We fix some $j$ and consider $\hat\Theta
_{\hat
\beta,j}$ as the $j$th row of
the nodewise regression $\hat{\Theta} = \hat{\Theta}_{\hat{\beta
}}$ in
(\ref{nodelassoglm}) based on the matrix input $\hat\Sigma_{\hat
\beta}$.

We let $\Sigma_{\beta} = \EE[\bx_{{\beta}}^T \bx_{{\beta}}/ n]$
and define
$\Theta:= \Theta_{\beta^0}:= \Sigma_{\beta^0}^{-1}$ (assumed to exist).
Let $s_j:= s_{\beta^0, j} $ be the number of off-diagonal zeros of the
$j$th row of $\Theta_{\beta^0} $.
Analogous to Section~\ref{subsecsub-Gaussian},\vspace*{1pt} we let $\bx_{\beta
^0, -j}
\gamma_{{\beta^0},j} $ be the
projection of $\bx_{{\beta^0}, j} $ on $\bx_{{\beta^0}, -j}$ using the
inner products in the matrix $\Sigma_{\beta^0}$ and let $\eta_{\beta^0,
j}:= \bx_{{\beta^0}, j }- \bx_{{\beta^0}, -j} \gamma_{\beta^0,j}$.
We then make the following assumptions:
\begin{longlist}[(D1)]
\item[(D1)] The pairs of random variables $\{(y_i, x_i)
\}_{i=1}^n $ are i.i.d. and $\|\bx\|_{\infty} = \max_{i,j}|\bx
_{i,j}| =
{\mathcal O}(K)$ and
$\|\bx_{\beta^0, -j} \gamma_{\beta^0, j} \|_{\infty}
= {\mathcal O} (K) $ for some $K \ge1$.

\item[(D2)] It holds that $K^2 s_j
\sqrt{\log(p)/ n } = o (1)$.

\item[(D3)] The smallest eigenvalue
of $\Sigma_{\beta^0} $ is bounded away from zero, and moreover, $\|
\Sigma_{\beta^0}
\|_{\infty}= {\mathcal O} (1) $.

\item[(D4)] For some $\delta>0$ and all $\| \beta- \beta^0 \|_1 \le
\delta$, it holds that
$w_{\beta} $ stays away from zero and that $\| w_{\beta} \|_{\infty} =
{\mathcal O} (1) $. We further require that for all such $\beta$ and
all $x$ and~$y$
\[
\bigl| w_{\beta}(y, x) - w_{\beta^0} (y,x) \bigr| \le\bigl| x \bigl( \beta-
\beta^0\bigr) \bigr|. %
\]

\item[(D5)]
It holds that
\[
\bigl\| \bx\bigl( \hat\beta- \beta^0\bigr) \bigr\|_n = {\mathcal
O}_{\PP} ( \lambda \sqrt{s_0 } ),\qquad\bigl \| \hat\beta-
\beta^0 \bigr\|_1 = {\mathcal O}_{\PP} ( \lambda
s_0 ). %
\]
\end{longlist}
Condition (D5) and (C2) typically hold when $\lambda\sqrt{s_0} = o(1)$
with tuning parameter $\lambda
\asymp\sqrt{\log(p)/n}$ since the compatibility condition is then
inherited from (D3) (see also Section~\ref{subsecassglm}).
We have the following result.
%
\begin{theo} \label{node2theorem} Assume the conditions \textup{(D1)--(D5)}.
Then, using $\lambda_j \asymp  K \sqrt{\log(p)/ n }$ for
the nodewise
lasso $\hat{\Theta}_{\hat{\beta}, j }$.
\begin{eqnarray*}
\| \hat\Theta_{\hat\beta, j} - \Theta_{\beta^0, j } \|_1 &=& {
\mathcal O}_{\PP} \bigl( K s_j \sqrt{\log(p)/ n } \bigr) + {
\mathcal O}_{\PP} \bigl( K^2 s_0 \bigl(\bigl(
\lambda^2 / \sqrt{\log(p)/ n } \bigr) \vee\lambda\bigr) \bigr),
\\
\| \hat\Theta_{\hat\beta, j} - \Theta_{\beta^0, j } \|_2& =& {
\mathcal O}_{\PP} \bigl( K \sqrt{s_j \log(p)/ n } \bigr) + {
\mathcal O}_{\PP} \bigl( K^2 \sqrt s_0 \lambda
\bigr), %
\end{eqnarray*}
and for $\tau_{\beta^0, j}^2:= \Theta_{\beta^0, j,j } $
\[
\bigl| \hat\tau_{\hat\beta, j}^2 - \tau_{\beta^0, j}^2 \bigr|
= {\mathcal O}_{\PP} \bigl( K \sqrt{s_j \log(p)/ n } \bigr)
+ {\mathcal O}_{\PP} \bigl(K^2 \sqrt{s_0}
\lambda \bigr). %
\]
Moreover,
\begin{eqnarray*}
&&\bigl| \hat{\Theta}_{\hat{\beta},j} \Sigma_{\beta^0} \hat{\Theta
}_{\hat{\beta
},j}^T - \Theta_{\beta^0,j,j}\bigr | %
\\
&&\qquad\le\| \Sigma_{\beta^0} \|_{\infty} \| \hat{\Theta}_{\hat{\beta},j} -
\Theta_{\beta^0, j} \|_1^2 \wedge
\Lambda_\mathrm{max}^2 \| \hat{\Theta}_{\hat{\beta},j} - \Theta
_{\beta^0,
j} \|_2^2 + 2 \bigl| \hat\tau_{\hat\beta, j}^2
- \tau_{\beta^0, j}^2 \bigr|, %
\end{eqnarray*}
where $\Lambda_\mathrm{max}^2$ is the maximal eigenvalue of $\Sigma
_{\beta
^0} $.
\end{theo}

A proof using ideas for establishing Theorem \ref{nodetheorem} is given
in the supplemental article \cite{supplement}.

\begin{corr} \label{node2corollary} Assume the conditions of Theorem
\ref{node2theorem},
with tuning parameter $\lambda\asymp\sqrt{\log(p)/ n}$, $K \asymp1$,
$s_j = o ( \sqrt n / \log(p)) $ and $s_0= o (\sqrt n / \log(p)) $.
Then
\begin{eqnarray*}
\| \hat\Theta_{\hat\beta, j} - \Theta_{\beta^0, j } \|_1 &=&
{o}_{\PP} \bigl( 1/ \sqrt{\log(p)} \bigr), %
\\
\| \hat\Theta_{\hat\beta, j} - \Theta_{\beta^0, j } \|_2 &=&
o_{\PP} \bigl(n^{-1/4}\bigr)
\end{eqnarray*}
and
\[
\bigl| \hat{\Theta}_{\hat{\beta},j} \Sigma_{\beta^0} \hat{\Theta
}_{\hat{\beta
},j}^T - \Theta_{\beta^0, j,j}\bigr | = o_{\PP}
\bigl( 1/ \log (p) \bigr). %
\]
\end{corr}

The next lemma is useful when estimating the asymptotic variance.

\begin{lemm} \label{normallemma} Assume the conditions of Corollary
\ref{node2corollary}.
Let for $i=1, \ldots, n$, $\xi_i$ be a real-valued random variable
and $x_i^T \in\R^p$, and let
$( x_i, \xi_i)_{i=1}^n$ be i.i.d. Assume $\EE x_i^T \xi_i=0 $ and
that $| \xi_i | \le1 $.
Then
\[
\hat{\Theta}_{\hat{\beta},j} \sum_{i=1}^n
x_i^T \xi_i / n = \Theta _{\beta^0,j}
\sum_{i=1}^n x_i^T
\xi_i /n + {o}_{\PP} \bigl(n^{-1/2} \bigr).
\]
Let $A:= \EE x_i^T x_i \xi_i^2$ (assumed to exist). Assume that
$\| A \Theta_j^T \|_{\infty} = {\mathcal O} (1) $ and that $1/
(\Theta
_j A \Theta_j^T) = {\mathcal O} (1)$.
Then
\[
\hat{\Theta}_{\hat{\beta},j} A \hat{\Theta}_{\hat{\beta},j}^T =
\Theta_{\beta^0,j} A \Theta_{\beta^0,j}^T + o_{\PP} (
1). %
\]
Moreover, then
\[
{ \hat{\Theta}_{\hat{\beta},j} \sum_{i=1}^n x_i^T \xi_i /\sqrt
{n} \over
\sqrt{\hat{\Theta}_{\hat{\beta},j} A \hat{\Theta}_{\hat{\beta
},j}^T } } %
\]
convergences weakly to a ${\cal N} (0,1 )$-distribution.
\end{lemm}

A proof is given in the supplemental article \cite{supplement}.

\subsubsection{Consequence for GLMs}\label{subsecassglm}
Consider the case where $ a \mapsto\rho(y, a)$ is convex for all $y$.
We let $\{ (y_i, x_i ) \}_{i=1}^n \sim\PP$ be i.i.d. random
variables. We denote by
$\bx_{\beta^0}$ the weighted design matrix\vadjust{\goodbreak}
$W_{\beta^0} \bx$ with $W_{\beta^0}$ the diagonal matrix with elements
$\{ \sqrt{\ddot\rho(y_i, x_i \beta^0 ) }
\}_{i=1}^n $. We further let
$\bx_{\beta^0, -j} \gamma_{\beta^0, j }^0$ be the projection in $L_2
(\PP)$
of $\mathbf{X}_{\beta^0, j} $ on $\bx_{\beta^0, -j}$, $j=1,
\ldots,
p $. We write
$\Sigma_{\beta^0}:= \EE{\bx}_{\beta^0}^T \mathbf{X}_{\beta^0} / n$
and let $s_j$ be the number
of nonzero lower-diagonal elements of the $j$th column of $\Sigma
_{\beta
^0}$ ($j=1, \ldots, p$).

\begin{theo} \label{GLMtheorem} Let $\{ (y_i, x_i ) \}_{i=1}^n \sim
\PP$ be i.i.d. random variables.
Assume:

\begin{longlist}[(vii)]
\item[(i)] Condition \textup{(C1)},
\item[(ii)] $\| 1 / \ddot\rho_{\beta^0} \|_{\infty} = {\mathcal
O } (1) $,
\item[(iii)] $\| \bx\|_{\infty} ={\mathcal O}(1)$,
\item[(iv)] $\|\bx\beta^0 \|_{\infty} = {\mathcal O} (1) $ and $\|
\bx
_{\beta^0, -j} \gamma_{\beta^0, j }^0 \|_{\infty} = {\mathcal O}
(1) $
for each $j$,
\item[(v)] the smallest eigenvalue of $\Sigma_{\beta^0} $ stays away
from zero,
\item[(vi)] $1/ (\Theta_{\beta^0,j} P \dot\rho_{\beta^0} \dot
\rho
_{\beta^0}^T\Theta_{\beta^0, j}^T) ={\mathcal O} (1)$ $\forall j$,
\item[(vii)] $s_0 = o( \sqrt n / \log(p)) $ and $s_j = \sqrt{ n /
\log(p)}$ for all $j$.
\end{longlist}
Take $\hat\Theta$ equal to $\hat\Theta_{\mathrm{Lasso}} $ given in
(\ref{nodelassoglm})
with $\lambda_j \asymp\sqrt{\log(p)/n}$ ($j=1, \ldots, p$)
suitably chosen.
For the estimator in (\ref{estglm}), with suitable $\lambda\asymp
\sqrt{\log(p) / n } $,
we have for each~$j$
\[
\sqrt{n}\bigl(\hat{b}_j - \beta_{j}^0\bigr)/
\hat{\sigma}_j = V_j + o_\mathbf {P} (1),
\]
where $V_j $ converges weakly to a ${\cal N}(0,1)$-distribution and
where
\[
\hat\sigma_j^2:= \bigl( \hat\Theta P_n
\dot\rho_{\hat\beta} \dot \rho_{\hat\beta}^T \hat{
\Theta}^T \bigr)_{j,j}. %
\]
\end{theo}

A proof is given in Section~\ref{GLMtheoremsection}.

Note that for the case where $\rho_{\beta}$ is the minus
log-likelihood, $P \dot\rho_{\beta^0} \dot\rho_{\beta^0}^T=
\Sigma
_{\beta^0}$,
and hence $ \Theta_{\beta^0,j} P \dot\rho_{\beta^0} \dot\rho
_{\beta
^0}^T\Theta_{\beta^0, j}^T =
\Theta_{\beta^0, j, j }$. Assumption {(vi)} then follows from
assumptions (i)--(iii) since
$1/ \Theta_{\beta^0, j, j } \le\Sigma_{\beta^0, j,j } $.

\section{Empirical results}\label{empiricalsection}

We consider finite sample behavior for inference of individual regression
coefficients $\beta^0_j$, including adjustment for the case of multiple
hypothesis testing.

\subsection{Methods and models}

We compare our method based on $\hat{b}_{\mathrm{Lasso}}$ with a procedure
based on multiple sample splitting \cite{memepb09} (for multiple
hypothesis testing
only) and with a residual bootstrap method proposed by \cite{chatter13}.

The implementational details for inference based on
$\hat{b}_{\mathrm{Lasso}}$ are as follows. For the linear regression
of the
response $Y$ versus the design $\bx$, we
use the scaled lasso~\cite{sunzhang11} with its universal
regularization parameter, and we use its estimate
$\hat{\sigma}^2_{\eps}$ of the error variance. For logistic
regression, we
use the corresponding lasso estimator with tuning parameter from 10-fold
cross-validation. Regarding the nodewise
lasso (for linear and logistic regression), we choose the same tuning
parameter $\lambda_j \equiv\lambda_X$ by
10-fold cross-validation among all nodewise regressions.
An alternative method which we did not yet examine in the simulations
would be to do nodewise regression with
square-root lasso using a universal choice for the tuning parameter
(see Remark~\ref{studentizeremark}).
For the bootstrap
method from \cite{chatter13}, we use 10-fold cross-validation to
sequentially select the tuning parameter for lasso and subsequently for
adaptive lasso. For multiple sample splitting \cite{memepb09}, we do
variable screening with the lasso whose regularization parameter is chosen
by 10-fold cross-validation.

The construction of confidence intervals and hypothesis tests for
individual parameters $\beta^0_j$ based on
$\hat{b}_{\mathrm{Lasso}}$ is straightforward, as described in
Section~\ref{subsecdespars}. Adjustment for multiple testing of hypotheses
$H_{0,j}$ over all $j=1,\ldots,p$ is done using the Bonferroni--Holm
procedure for controlling the family-wise error rate (FWER). For the
bootstrap procedure from \cite{chatter13}, the Bonferroni--Holm adjustment
is not sensible, unless we would draw very many bootstrap resamples
(e.g., 10,000
or more): with fewer resamples, we cannot reliably estimate the
distribution in the tails needed for Bonferroni--Holm correction. Thus, for
this bootstrap method, we
only consider construction of confidence intervals. Finally, the
multiple sample splitting method \cite{memepb09} is
directly giving $p$-values which control the FWER.

For our simulation study, we consider (logistic) linear
models where the rows of
$\bx$ are fixed i.i.d. realizations from ${\cal N}_p(0,\Sigma)$. We specify
two different covariance matrices:
\begin{eqnarray*}
&&\mbox{Toeplitz:}\qquad \Sigma_{j,k} = 0.9^{|j-k|},
\\
&&\mbox{Equi corr:}\qquad\Sigma_{j,k}
\equiv0.8\qquad \mbox{for all } j \neq k,\qquad \Sigma_{j,j} \equiv1\qquad
\mbox{for all }
j.
\end{eqnarray*}
%
The active set has either cardinality $s_0 = |S_0| = 3$ or $s_0 = 15$, and
each of it is of one of the following forms:
\[
S_0 = \{1,2,\ldots,s_0\}, \mbox{or: realization of
random support } S_0= \{u_1,\ldots,u_{s_0}\},
\]
where $u_1,\ldots,u_{s_0}$ is a fixed realization of $s_0$ draws without
replacement from $\{1,\ldots,p\}$. The regression coefficients are from
a fixed realization of $s_0$ i.i.d. Uniform $U[0,c]$ variables with $c
\in
\{1,2,4\}$.
For linear models, the distribution of the errors is always $\eps
_1,\ldots,\eps_n \sim{\cal N}(0,1)$; see comment below regarding
$t$-distributed errors. We also consider logistic regression models with
binary response and
\[
\log\bigl(\pi(x)/\bigl(1 - \pi(x)\bigr)\bigr) = x \beta^0,\qquad \pi(x) =
\PP[y_1=1|x_1=x].
\]
Sample size is always $n = 100$ (with some exceptions in the
supplemental article~\cite{supplement}) and the number of variables is
$p = 500$. We
then consider many combinations of the different specifications above. All
our results are based on 100 independent simulations of the model with
fixed design and fixed regression
coefficients (i.e.,
repeating over 100 independent simulations of the errors in a linear model).

\subsection{Results for simulated data}

\subsubsection{Linear model: Confidence intervals}
We consider average coverage and average length of
the intervals for individual coefficients corresponding to variables in
either $S_0$ or $S_0^c$: denoting by $\mathrm{CI}_j$ a two-sided
confidence interval for $\beta^0_j$, we report empirical versions of
\begin{eqnarray*}
\operatorname{Avgcov} S_0 &=& s_0^{-1}
\sum_{j \in S_0} \PP \bigl[\beta_j^0
\in \mathrm{CI}_j\bigr], \\
\operatorname{Avgcov} S_0^c
&=& (p-s_0)^{-1} \sum_{j \in S_0^c}
\PP[0 \in \mathrm{CI}_j],
\\
\operatorname{Avglength} S_0 &=& s_0^{-1}
\sum_{j \in S_0} \mathrm{length}(\mathrm{CI}_j);\qquad
\mbox{and analogously for } \operatorname {Avglength} S_0^c.
\end{eqnarray*}
The following Tables~\ref{tab1}--\ref{tab4} are for different active sets.

%
\begin{table}
\caption{Linear model: average coverage and length of confidence
intervals, for nominal coverage equal to $0.95$. ``Lasso-Pro''
(lasso-projection) denotes the procedure based on our desparsified
estimator $\hat{b}_{\mathrm{Lasso}}$; ``Res-Boot'' is the residual based
bootstrap from \protect\cite{chatter13}}
\label{tab1}
%
\begin{tabular*}{\textwidth}{@{\extracolsep{\fill}}lcd{1.3}d{1.3}d{1.3}d{1.3}@{}}
\hline
& & \multicolumn{2}{c}{\textbf{Toeplitz}}  & \multicolumn{2}{c@{}}{\textbf{Equi corr}} \\[-6pt]
& & \multicolumn{2}{c}{\hrulefill}  & \multicolumn{2}{c@{}}{\hrulefill} \\
\multicolumn{1}{@{}l}{\textbf{Measure}} &
\multicolumn{1}{c}{\textbf{Method}} & \multicolumn{1}{c}{$\bolds{U([0,2])}$} & \multicolumn{1}{c}{$\bolds{U([0,4])}$} &
\multicolumn{1}{c}{$\bolds{U([0,2])}$} & \multicolumn{1}{c@{}}{$\bolds{U([0,4])}$} \\
\hline
\multicolumn{6}{c}{Active set $S_0 = \{1,2,3\}$}\\
Avgcov $S_0$ & Lasso-Pro & 0.86 & 0.84 & 0.90 & 0.89 \\
& Res-Boot & 0.66 & 0.85 & 0.45 & 0.57 \\
Avglength $S_0$ & Lasso-Pro & 0.786 & 0.787 & 0.762 & 0.760 \\
& Res-Boot & 0.698 & 0.918 & 0.498 & 0.670 \\[3pt]
Avgcov $S_0^c$ & Lasso-Pro & 0.95 & 0.95 & 0.95 & 0.95 \\
& Res-Boot & 1.00 & 1.00 & 1.00 & 1.00 \\
Avglength $S_0^c$ & Lasso-Pro & 0.786 & 0.787 & 0.811 & 0.808 \\
& Res-Boot & 0.000 & 0.000 & 0.006 & 0.007 \\
\hline
\end{tabular*}
\end{table}

\begin{table}
\caption{See caption of Table \protect\ref{tab1}}\label{tab2}
\begin{tabular*}{\textwidth}{@{\extracolsep{\fill}}lcd{1.3}d{1.3}d{1.3}d{1.3}@{}}
\hline
& & \multicolumn{2}{c}{\textbf{Toeplitz}}  & \multicolumn{2}{c@{}}{\textbf{Equi corr}} \\[-6pt]
& & \multicolumn{2}{c}{\hrulefill}  & \multicolumn{2}{c@{}}{\hrulefill} \\
\multicolumn{1}{@{}l}{\textbf{Measure}} &
\multicolumn{1}{c}{\textbf{Method}} & \multicolumn{1}{c}{$\bolds{U([0,2])}$} & \multicolumn{1}{c}{$\bolds{U([0,4])}$} &
\multicolumn{1}{c}{$\bolds{U([0,2])}$} & \multicolumn{1}{c@{}}{$\bolds{U([0,4])}$} \\
\hline
\multicolumn{6}{c}{Active set with $s_0 = 3$ and support from fixed random realization}\\
Avgcov $S_0$ & Lasso-Pro & 0.95 & 0.94 & 0.89 & 0.87 \\
& Res-Boot & 0.58 & 0.73 & 0.31 & 0.51 \\
Avglength $S_0$ & Lasso-Pro & 0.890 & 0.934 & 0.822 & 0.821 \\
& Res-Boot & 0.336 & 0.463 & 0.500 & 0.743 \\[3pt]
Avgcov $S_0^c$ & Lasso-Pro & 0.95 & 0.96 & 0.95 & 0.95 \\
& Res-Boot & 1.00 & 1.00 & 1.00 & 1.00 \\
Avglength $S_0^c$ & Lasso-Pro & 0.879 & 0.923 & 0.805 & 0.804 \\
& Res-Boot & 0.002 & 0.003 & 0.008 & 0.009 \\
\hline
\end{tabular*}  \vspace*{-3pt}
\end{table}

\begin{table}
\caption{See caption of Table \protect\ref{tab1}}\label{tab3}
\begin{tabular*}{\textwidth}{@{\extracolsep{\fill}}lcd{1.3}d{1.3}d{1.3}d{1.3}@{}}
\hline
& & \multicolumn{2}{c}{\textbf{Toeplitz}}  & \multicolumn{2}{c@{}}{\textbf{Equi corr}} \\[-6pt]
& & \multicolumn{2}{c}{\hrulefill}  & \multicolumn{2}{c@{}}{\hrulefill} \\
\multicolumn{1}{@{}l}{\textbf{Measure}} &
\multicolumn{1}{c}{\textbf{Method}} & \multicolumn{1}{c}{$\bolds{U([0,2])}$} & \multicolumn{1}{c}{$\bolds{U([0,4])}$} &
\multicolumn{1}{c}{$\bolds{U([0,2])}$} & \multicolumn{1}{c@{}}{$\bolds{U([0,4])}$} \\
\hline
\multicolumn{6}{c}{Active set $S_0 = \{1,2,\ldots,15\}$}\\
Avgcov $S_0$ & Lasso-Pro & 0.76 & 0.73 & 0.56 & 0.53 \\
& Res-Boot & 0.79 & 0.87 & 0.54 & 0.63 \\
Avglength $S_0$ & Lasso-Pro & 0.813 & 0.814 & 0.559 & 0.554 \\
& Res-Boot & 1.012 & 1.138 & 0.746 & 0.844 \\[3pt]
Avgcov $S_0^c$ & Lasso-Pro & 0.96 & 0.96 & 0.93 & 0.93 \\
& Res-Boot & 1.00 & 1.00 & 0.99 & 0.99 \\
Avglength $S_0^c$ & Lasso-Pro & 0.788 & 0.789 & 0.568 & 0.562 \\
& Res-Boot & 0.001 & 0.000 & 0.047 & 0.044 \\
\hline
\end{tabular*}      \vspace*{-3pt}
\end{table}

\begin{table}8
\caption{See caption of Table \protect\ref{tab1}}\label{tab4}
%
\begin{tabular*}{\textwidth}{@{\extracolsep{\fill}}lcd{1.3}d{1.3}d{1.3}d{1.3}@{}}
\hline
& & \multicolumn{2}{c}{\textbf{Toeplitz}}  & \multicolumn{2}{c@{}}{\textbf{Equi corr}} \\[-6pt]
& & \multicolumn{2}{c}{\hrulefill}  & \multicolumn{2}{c@{}}{\hrulefill} \\
\multicolumn{1}{@{}l}{\textbf{Measure}} &
\multicolumn{1}{c}{\textbf{Method}} & \multicolumn{1}{c}{$\bolds{U([0,2])}$} & \multicolumn{1}{c}{$\bolds{U([0,4])}$} &
\multicolumn{1}{c}{$\bolds{U([0,2])}$} & \multicolumn{1}{c@{}}{$\bolds{U([0,4])}$} \\
\hline
\multicolumn{6}{c}{Active set with $s_0 = 15$ and support from fixed random realization}\\
Avgcov $S_0$ & Lasso-Pro & 0.93 & 0.94 & 0.55 & 0.44 \\
& Res-Boot & 0.45 & 0.54 & 0.48 & 0.55 \\
Avglength $S_0$ & Lasso-Pro & 2.391 & 4.354 & 0.572 & 0.552 \\
& Res-Boot & 0.480 & 0.599 & 0.675 & 0.809 \\[3pt]
Avgcov $S_0^c$ & Lasso-Pro & 0.95 & 0.95 & 0.93 & 0.92 \\
& Res-Boot & 0.98 & 0.97 & 0.99 & 0.99 \\
Avglength $S_0^c$ & Lasso-Pro & 2.370 & 4.317 & 0.570 & 0.550 \\
& Res-Boot & 0.029 & 0.035 & 0.048 & 0.050 \\
\hline
\end{tabular*}
\end{table}

\emph{Discussion}.
As the main finding, we summarize that the desparsified lasso estimator is clearly
better for the variables in $S_0$ than the residual based bootstrap. For
the variables in $S_0^c$ with regression coefficients equal to zero, the
residual bootstrap exhibits the super-efficiency phenomenon: the average
length\vadjust{\goodbreak} of the interval is often very close to zero while coverage equals
one. This cannot happen with the desparsified lasso estimator: in
contrast to
the residual based bootstrap, the desparsified lasso estimator allows
for a
convergence result which is uniform for a large class of parameters, and
hence leading to honest confidence intervals; see Section~\ref
{uniformsection}. Furthermore, our empirical results for active sets
with $s_0 = 15$ indicate that inference with the desparsified lasso has its
limit when the problem is not sufficiently sparse, especially for the case
with equi-correlated design: this is in line with our theoretical results.

%
\begin{table}
\caption{Linear model: family-wise error rate (FWER) and power of multiple
testing, for nominal FWER equal to
$0.05$. ``Lasso-Pro''(lasso-projection) denotes the procedure based on
our de-sparsified estimator $\hat{b}_{\mathrm{Lasso}}$ with
Bonferroni--Holm adjustment for multiple testing; ``MS-Split'' is
the multiple sample splitting method from \cite{memepb09}}\label{tab5}
%
\begin{tabular*}{\textwidth}{@{\extracolsep{\fill}}lccccc@{}}
\hline
& & \multicolumn{2}{c}{\textbf{Toeplitz}}  & \multicolumn{2}{c@{}}{\textbf{Equi corr}} \\[-6pt]
& & \multicolumn{2}{c}{\hrulefill}  & \multicolumn{2}{c@{}}{\hrulefill} \\
\multicolumn{1}{@{}l}{\textbf{Measure}} &
\multicolumn{1}{c}{\textbf{Method}} & \multicolumn{1}{c}{$\bolds{U([0,2])}$} & \multicolumn{1}{c}{$\bolds{U([0,4])}$} &
\multicolumn{1}{c}{$\bolds{U([0,2])}$} & \multicolumn{1}{c@{}}{$\bolds{U([0,4])}$} \\
\hline
\multicolumn{6}{c}{Active set $S_0 = \{1,2,3\}$}\\
Power & Lasso-Pro & 0.42 & 0.69 & 0.48 & 0.82 \\
& MS-Split & 0.60 & 0.83 & 0.35 & 0.63 \\
FWER & Lasso-Pro & 0.03 & 0.05 & 0.13 & 0.13 \\
& MS-Split & 0.16 & 0.25 & 0.00 & 0.00 \\
\hline
\end{tabular*}   \vspace*{-3pt}
\end{table}

\begin{table}[b]\vspace*{-3pt}
\caption{See caption of Table~\protect\ref{tab5}}\label{tab6}
%
\begin{tabular*}{\textwidth}{@{\extracolsep{\fill}}lccccc@{}}
\hline
& & \multicolumn{2}{c}{\textbf{Toeplitz}}  & \multicolumn{2}{c@{}}{\textbf{Equi corr}} \\[-6pt]
& & \multicolumn{2}{c}{\hrulefill}  & \multicolumn{2}{c@{}}{\hrulefill} \\
\multicolumn{1}{@{}l}{\textbf{Measure}} &
\multicolumn{1}{c}{\textbf{Method}} & \multicolumn{1}{c}{$\bolds{U([0,2])}$} & \multicolumn{1}{c}{$\bolds{U([0,4])}$} &
\multicolumn{1}{c}{$\bolds{U([0,2])}$} & \multicolumn{1}{c@{}}{$\bolds{U([0,4])}$} \\
\hline
\multicolumn{6}{c}{Active set with $s_0=3$ and support from fixed random realization}\\
Power & Lasso-Pro & 0.54 & 0.81 & 0.56 & 0.79 \\
& MS-Split & 0.44 & 0.71 & 0.40 & 0.69 \\
FWER & Lasso-Pro & 0.00 & 0.00 & 0.10 & 0.11 \\
& MS-Split & 0.00 & 0.00 & 0.00 & 0.00 \\
\hline
\end{tabular*}
\end{table}

Finally, we have also looked at non-Gaussian models where the error terms
are from a scaled $t_5$ distribution (Student distribution with 5
degrees of freedom)
with variance equal to
one. The results (not reported here) look essentially
identical as in Tables~\ref{tab1}--\ref{tab4}.

\subsubsection{Linear model: Multiple testing}
We consider multiple two-sided testing of hypotheses $H_{0,j}; \beta
^0_j =
0$ among all $j=1,\ldots,p$. We correct the $p$-values based on our
$\hat{b}_{\mathrm{Lasso}}$ with the Bonferroni--Holm procedure to control
the familywise error rate (FWER). The method based on multiple sample splitting
\cite{memepb09} automatically yields $p$-values for controlling the
FWER. For measuring power, we report on the empirical version of
\[
\mathrm{Power} = s_0^{-1} \sum
_{j \in S_0} \PP[H_{0,j} \mbox{ is rejected}].
\]
The following Tables~\ref{tab5}--\ref{tab8} are for different active sets.

\begin{table}
\caption{See caption of Table~\protect\ref{tab5}}\label{tab7}
%
\begin{tabular*}{\textwidth}{@{\extracolsep{\fill}}lccccc@{}}
\hline
& & \multicolumn{2}{c}{\textbf{Toeplitz}}  & \multicolumn{2}{c@{}}{\textbf{Equi corr}} \\[-6pt]
& & \multicolumn{2}{c}{\hrulefill}  & \multicolumn{2}{c@{}}{\hrulefill} \\
\multicolumn{1}{@{}l}{\textbf{Measure}} &
\multicolumn{1}{c}{\textbf{Method}} & \multicolumn{1}{c}{$\bolds{U([0,2])}$} & \multicolumn{1}{c}{$\bolds{U([0,4])}$} &
\multicolumn{1}{c}{$\bolds{U([0,2])}$} & \multicolumn{1}{c@{}}{$\bolds{U([0,4])}$} \\
\hline
\multicolumn{6}{c}{Active set $S_0 = \{1,2,\ldots, 15\}$}\\
Power & Lasso-Pro & 0.73 & 0.89 & 0.70 & 0.92 \\
& MS-Split & 0.23 & 0.67 & 0.00 & 0.00 \\
FWER & Lasso-Pro & 0.03 & 0.02 & 1.00 & 1.00 \\
& MS-Split & 0.00 & 0.00 & 0.00 & 0.00 \\
\hline
\end{tabular*}\vspace*{-4pt}
\end{table}

\begin{table}
\caption{See caption of Table~\protect\ref{tab5}}\label{tab8}
%
\begin{tabular*}{\textwidth}{@{\extracolsep{\fill}}lccccc@{}}
\hline
& & \multicolumn{2}{c}{\textbf{Toeplitz}}  & \multicolumn{2}{c@{}}{\textbf{Equi corr}} \\[-6pt]
& & \multicolumn{2}{c}{\hrulefill}  & \multicolumn{2}{c@{}}{\hrulefill} \\
\multicolumn{1}{@{}l}{\textbf{Measure}} &
\multicolumn{1}{c}{\textbf{Method}} & \multicolumn{1}{c}{$\bolds{U([0,2])}$} & \multicolumn{1}{c}{$\bolds{U([0,4])}$} &
\multicolumn{1}{c}{$\bolds{U([0,2])}$} & \multicolumn{1}{c@{}}{$\bolds{U([0,4])}$} \\
\hline
\multicolumn{6}{c}{Active set with $s_0=15$ and support from fixed random realization}\\
Power & Lasso-Pro & 0.06 & 0.07 & 0.65 & 0.86 \\
& MS-Split & 0.07 & 0.14 & 0.00 & 0.00 \\
FWER & Lasso-Pro & 0.02 & 0.00 & 0.96 & 0.98 \\
& MS-Split & 0.02 & 0.13 & 0.00 & 0.00 \\
\hline
\end{tabular*}   \vspace*{-4pt}
\end{table}

\emph{Discussion.} Similarly to what we found for confidence intervals
above, multiple testing with the desparsified lasso estimator is
reliable and
works well for sparse problems (i.e., $s_0 = 3$). For less
sparse problems (i.e., $s_0 = 15$), the error control is less reliable,
especially for equi-correlated designs. For sparse Toeplitz designs,
the lasso-projection method has more power than multiple sample splitting,
a finding which is in line with our established optimality theory.

\subsubsection{Logistic regression: Multiple
testing}\label{subsecglmmulttest}\vspace*{-2pt}

The residual bootstrap\break  method~\cite{chatter13}
cannot be used in a straightforward way for logistic regression. As for
linear models, we compare our desparsified lasso estimator with the
multiple sample
splitting procedure, in the context of multiple testing for controlling the
FWER.

For the case of logistic regression shown in Tables~\ref{tab9}--\ref{tab10},
inference with the de-sparsified lasso method is not very reliable with
respect to the FWER. The multiple sample splitting method is found to
perform better. We present in the supplemental article \cite
{supplement} some additional
results for sample sizes $n=200$ and $n=400$, illustrating that the FWER
control as well as the power for the desparsified lasso
improve.\vspace*{-2pt}

\subsection{Real data analysis}

We consider a dataset about
riboflavin (vitamin~$B_{2}$) production\vadjust{\goodbreak} by bacillus subtilis. The data
has been
kindly provided by DSM (Switzerland) and is publicly available \cite
{bumeka13}.
The real-valued response variable is the logarithm of the riboflavin
production rate and there are $p=4088$
covariates (genes) measuring the logarithm of the expression level of 4088
genes. These measurements are from $n=71$ samples of genetically
engineered mutants of bacillus subtilis. We model the data with a
high-dimensional linear model and obtain the following results for
significance. The desparsified lasso procedure finds no significant
coefficient while the multiple sample splitting method claims
significance of one variable at the 5\% significance level for the
FWER. Such low power is to be expected in presence of thousands of
variables: finding significant groups of highly correlated variables would
seem substantially easier, at the price of not being able to infer
significant of variables at the individual level.

%
%

%

%
\begin{table}
\caption{Logistic regression: All other specifications as in
Table~\protect\ref{tab5}}\label{tab9}
\begin{tabular*}{\textwidth}{@{\extracolsep{\fill}}lcccc@{}}
\hline
& & \multicolumn{3}{c@{}}{\textbf{Toeplitz}}\\[-6pt]
& & \multicolumn{3}{c@{}}{\hrulefill} \\
\multicolumn{1}{@{}l}{\textbf{Measure}} &
\multicolumn{1}{c}{\textbf{Method}} & \multicolumn{1}{c}{$\bolds{U([0,1])}$} & \multicolumn{1}{c}{$\bolds{U([0,2])}$} &
\multicolumn{1}{c@{}}{$\bolds{U([0,4])}$}\\
\hline
%
Power & Lasso-ProG & 0.06 & 0.27 & 0.50 \\
& MS-Split & 0.07 & 0.37 & 0.08 \\[3pt]
FWER & Lasso-ProG & 0.03 & 0.08 & 0.23 \\
& MS-Split & 0.01 & 0.00 & 0.00 \\
\hline
\end{tabular*}   \vspace*{-4pt}
\end{table}

\begin{table}
\caption{Logistic regression: All other specifications as in
Table~\protect\ref{tab5}}\label{tab10}
%
\begin{tabular*}{\textwidth}{@{\extracolsep{\fill}}lcccc@{}}
\hline
& & \multicolumn{3}{c@{}}{\textbf{Toeplitz}}\\[-6pt]
& & \multicolumn{3}{c@{}}{\hrulefill} \\
\multicolumn{1}{@{}l}{\textbf{Measure}} &
\multicolumn{1}{c}{\textbf{Method}} & \multicolumn{1}{c}{$\bolds{U([0,1])}$} & \multicolumn{1}{c}{$\bolds{U([0,2])}$} &
\multicolumn{1}{c@{}}{$\bolds{U([0,4])}$}\\
\hline
Power & Lasso-ProG & 0.02 & 0.16 & 0.35 \\
& MS-Split & 0.00 & 0.17 & 0.27 \\[3pt]
FWER & Lasso-ProG & 0.08 & 0.16 & 0.27 \\
& MS-Split & 0.00 & 0.03 & 0.01 \\
\hline
\end{tabular*}  \vspace*{-4pt}
\end{table}
%

\section{Proofs and materials needed}\label{secproofs}

\subsection{Bounds for \texorpdfstring{$\|\hat{\beta}-\beta^0\|_1$}{||beta-beta0||1} with fixed
design}\label{subsecell1bound}

The following known result gives a bound for the $\ell_1$-norm estimation
accuracy.
%
\begin{lemm}\label{lem1}
Assume a linear model as in (\ref{modlin}) with Gaussian error and fixed
design $\bx$ which satisfies the compatibility condition with compatibility
constant $\phi_0^2$ and with $\hat{\Sigma}_{j,j} \le M^2 < \infty$
for all
$j$. Consider the lasso with regularization parameter $\lambda\ge2 M
\sigma_{\eps} \sqrt{\frac{2(t^2 + \log(p))}{n}}$.
Then, with probability at least $1 - 2 \exp(-t^2)$,
\[
\bigl\|\hat{\beta} - \beta^0\bigr\|_1 \le8 \lambda
\frac{s_0}{\phi_0^2} \quad\mbox {and}\quad\bigl \|\bx\bigl(\hat{\beta} - \beta^0\bigr)
\bigr\|_2^2/n \le8 \lambda^2 \frac
{s_0}{\phi_0^2}.
\]
\end{lemm}
A proof follows directly from the arguments in \cite{pbvdg11}, Theorem~6.1, which
can be modified to treat the case with unequal values of
$\hat{\Sigma}_{j,j}$ for various $j$. 

\subsection{Proof of Theorem \protect\ref{th1}}\label{subsecproofth1}

It is straightforward to see that
%
\begin{eqnarray}
\label{prlabel1} \|\Delta\|_{\infty} /\sqrt n &=&\bigl \|(\hat{
\Theta}_{\mathrm{Lasso}} \hat {\Sigma} - I) \bigl(\hat{\beta} - \beta^0
\bigr)\bigr\|_{\infty}
\nonumber
\\[-8pt]
\\[-8pt]
\nonumber
&\le&\bigl\|(\hat{\Theta }_{\mathrm{Lasso}} \hat{\Sigma} - I)
\bigr\|_{\infty} \bigl\|\hat{\beta} - \beta^0\bigr\|_1.
\end{eqnarray}
Therefore, by (\ref{add-inver}) we have that
$\|\Delta\|_{\infty} \le\sqrt n \| \hat\beta- \beta^0 \|_1 \max_j
{\lambda_j}/{\hat{\tau}_j^2} $, and using the bound from
Lemma \ref{lem1} completes the proof. 


\subsection{Random design: Bounds for compatibility constant and
\texorpdfstring{$\|\hat T^{-2}\|_{\infty}$}{||T-2||infinity}}\label{subseccompbound}

The compatibility condition with constant $\phi_0^2$ being bounded away
from zero
is ensured by a rather natural condition about sparsity.
We have the following result.
%
\begin{lemm}\label{lem2}
Assume \textup{(A2)}. Furthermore, assume that
$s_0 = o(n/\log(p))$.
Then there is a constant $L= {\mathcal O} (1) $ depending on $\Lambda
_\mathrm{min}$ only such that with probability tending to one the compatibility
condition holds with compatibility constant $\phi_{0}^2 \ge1/L^2 $.
\end{lemm}
A proof follows directly as in \cite{raskutti2010restricted}, Theorem 1.

Lemmas \ref{lem1} and \ref{lem2} say that we have a bound
%
\begin{eqnarray}
\label{ell1-bound} \bigl\|\hat{\beta} - \beta^0\bigr\|_1 &=& {\mathcal
O}_{\PP} \biggl(s_0 \sqrt {\frac
{\log(p)}{n}} \biggr),
\nonumber
\\[-8pt]
\\[-8pt]
\nonumber
\bigl\|
\bx\bigl(\hat{\beta} -\beta^0\bigr)\bigr\|_2^2/n& =&{
\mathcal O}_{\PP} \biggl( \frac{s_0 \log(p)}{n} \biggr),
\end{eqnarray}
when assuming (A2) and sparsity $s_0 =
o(n/\log(p))$.

When using the lasso for nodewise regression in (\ref{thetalasso}), we would
like to have a bound for $\|\hat{T}_{\mathrm{Lasso}}^{-2}\|_{\infty}$
appearing in Theorem \ref{th1}.
%
\begin{lemm}\label{lem3}
Assume \textup{(A2)} with row-sparsity for $\Theta:= \Sigma^{-1}$ bounded by
\[
\max_j s_j = o\bigl(n / \log(p)\bigr).
\]
%
Then, when suitably choosing the regularization parameters $\lambda_j
\asymp\sqrt{\log(p)/n}$
uniformly in $j$,
%
\[
\max_j 1/ \hat\tau_j^2 = {
\mathcal O}_{\PP}(1).
\]
\end{lemm}

\begin{pf}
A proof follows using standard arguments. With probability tending to
one the compatibility assumption holds
uniformly for all nodewise regressions with compatibility constant
bounded away from zero uniformly in $j$, as in Lemma~\ref{lem2} and
invoking the union bound.
Furthermore, the
population error variance $\tau_j^2 = \EE[(X_{1,j} - \sum_{k \neq j}
\gamma_{j,k} X_{1,k})^2]$, where $\gamma_{j,k}$ are the population regression
coefficients of $X_{1,j}$ versus $\{X_{1,k}; k \neq j\}$ satisfy:
uniformly in $j$,
$\tau_{j}^2 = 1/ \Theta_{j,j}\ge\Lambda_{\mathrm{min}}^2 > 0$
and $\tau_j^2 \le\EE[X_{1,j}^2] =
\Sigma_{j,j} = {\mathcal O}(1)$, thereby invoking
assumption (A2).
Thus, all the error
variances behave nicely.
Recall that
\[
\hat\tau_j^2:= \| X_j - {
\bx}_{-j} \hat\gamma_j \|_2^2 / n +
\lambda _j \| \hat\gamma_j \|_1.
\]
In the following, the probability statements are again uniformly in $j$
by the union bound for
suitable tuning parameters $\lambda_j \asymp\sqrt{\log(p) /n }$
uniformly in $j$.\vadjust{\goodbreak}
Each
nodewise regression satisfies $\|\bx_{-j}( \hat{\gamma}_j - \gamma
_j )\|
_2^2/n = {\mathcal O}_{\PP}(s_j
\log(p)/n)$ [see Lemma \ref{lem1} or (\ref{ell1-bound}), now applied
to the
lasso estimator for the regression of $X_j$ on $\mathbf{X}_{-j}$].
It follows that
\begin{eqnarray*}
\| X_j - {\bx}_{-j} \hat\gamma_j
\|_2^2 / n &= &\| X_j - {\bx}_{-j}
\gamma_j \|_2^2 / n + \bigl\|\bx_{-j}(
\hat{\gamma}_j - \gamma_j )\bigr\|_2^2/n
\\
&&{}+ 2 (X_j - {\bx}_{-j} \gamma_j
)^T \bx_{-j}( \hat{\gamma}_j -
\gamma_j ) %
\\
&=& \tau_j^2 + {\mathcal O}_{\PP} \bigl(
n^{-1/2}\bigr) + {\mathcal O}_{\PP} \bigl( \lambda_j^2
s_j \bigr) + {\mathcal O}_{\PP} ( \lambda_j
\sqrt{s_j} ) = \tau_j^2 + o_{\PP}
(1). %
\end{eqnarray*}

Note further that
\[
\| \gamma_j \|_1 \le\sqrt{s_j} \|
\gamma_j \|_2 \le\sqrt{s_j \Sigma
_{j,j}} / \Lambda_\mathrm{min}. %
\]
Moreover, by the same arguments giving the bounds in (\ref
{ell1-bound}), $\| \hat\gamma_j - \gamma_j \|_1 = {\mathcal O}_{\PP}
( s_j \lambda_j ) $ so that
\[
\lambda_j \| \hat\gamma_j \|_1 \le
\lambda_j \| \gamma_j \|_1 + \lambda
_j \| \hat\gamma_j - \gamma_j
\|_1 = \lambda_j {\mathcal O} (\sqrt {s_j}
) + \lambda_j {\mathcal O}_{\PP} ( \lambda_j
s_j ) = o_{\PP} (1).
\]
Hence, the statement of the lemma follows.
\end{pf} 

\subsection{Bounds for \texorpdfstring{$\|\hat{\beta}-\beta^0\|_2$}{||beta-beta0||2} with random
design}\label{subsecell2bound}

Note that $\| {\bx} (\hat\beta- \beta^0 ) \|_2^2 /\break  n = (\hat\beta-
\beta^0 )^T \hat\Sigma(\hat\beta- \beta^0 ) $.
Lemma \ref{lem2} uses \cite{raskutti2010restricted}, Theorem 1. The
same result can be invoked
to conclude that when (A2) holds and when $\lambda\asymp\sqrt{\log
(p) / n }$ is suitably chosen, then
for a suitably chosen fixed $C$,
with probability tending to
one
\begin{eqnarray*}
&&\bigl(\hat\beta- \beta^0 \bigr)^T \Sigma\bigl(\hat\beta-
\beta^0 \bigr) %
\\
&&\qquad\le\bigl(\hat\beta- \beta^0 \bigr)^T \hat\Sigma\bigl(
\hat\beta- \beta^0 \bigr) C + \sqrt{\frac{\log(p)} {n}}\bigl \| \hat\beta-
\beta^0\bigr \|_1 C.
\end{eqnarray*}
Hence,
\[
\bigl(\hat\beta- \beta^0 \bigr)^T \Sigma\bigl(\hat\beta-
\beta^0 \bigr) = {\mathcal O}_{\PP} \biggl(
\frac{s_0 \log(p) }{n} \biggr). %
\]
So under (A2) for suitable $\lambda\asymp\sqrt{\log(p) / n }$
%
\begin{equation}
\label{ell2-bound} \bigl\|\hat{\beta} - \beta^0\bigr\|_2 = {\mathcal
O}_{\PP}\bigl(\sqrt{s_{0} \log(p)/n}\bigr)
\end{equation}
(see also \cite{brt09}). This result will be applied in the next
subsection, albeit to the
lasso for node wise regression instead of for the original linear model.

\subsection{Proof of Theorem \protect\ref{th2}}\label{subsecproofth2}

Invoking Theorem \ref{th1} and Lemma \ref{lem3}, we have that
\[
\|\Delta\|_{\infty} \le{\mathcal O}_{\PP}\bigl(s_{0}
\log(p) / \sqrt n \bigr) = o_{\PP}(1),
\]
where the last bound follows by the sparsity assumption on $s_{0}$.

What remains to be shown is that $\|\hat\Omega- \Theta\|_{\infty} =
o_{\PP}(1)$, as detailed by the following lemma.
%
\begin{lemm}\label{lem4} Let $\hat\Theta:= \hat\Theta_{\mathrm
{Lasso}}$ with suitable tuning
parameters $\lambda_j$ satisfying $\lambda_j \asymp\sqrt{ \log(p)/
n} $ uniformly in $j$.
Assume the conditions of Lemma \ref{lem3}.
Suppose that $\max_j \lambda_j^2
s_j= o(1)$. Then
\[
\|\hat\Omega- \Theta\|_{\infty} = o_{\PP}(1).
\]
\end{lemm}


\begin{pf}
By the same arguments as in the proof of Lemma \ref{lem3}, uniformly
in $j$,
\[
\| \hat\Theta_j \|_1 = O_{\PP} (
\sqrt{s_j}).
\]
Furthermore, we have
%
\begin{equation}
\label{addprth2-2} \hat\Omega= \hat{\Theta} \hat{\Sigma} \hat{\Theta}^T =
(\hat {\Theta} \hat{\Sigma} - I) \hat{\Theta}^T + \hat{
\Theta}^T
\end{equation}
and
%
\begin{equation}
\label{addprth2-3}  \bigl\|(\hat{\Theta} \hat{\Sigma} - I) \hat{\Theta
}^T\bigr\| _{\infty} \le \max_j
\lambda_j \|\hat{\Theta}_j\|_1 / \hat
\tau_j^2 = o_{\PP}(1),
\end{equation}
which follows from Lemma \ref{lem3}.
Finally, we have using standard arguments for the $\ell_2$-norm bounds [see
also (\ref{ell2-bound})]
%
\begin{equation}
\label{addprth2-4} \|\hat{\Theta} - \Theta\|_{\infty} \le\max
_j \|\hat{\Theta}_j - \Theta_j
\|_2 \le\max_j \lambda_j
\sqrt{s_j} = o_{\PP}(1).
\end{equation}
Using (\ref{addprth2-2})--(\ref{addprth2-4}), we complete the
proof.
\end{pf}

The proof of Theorem \ref{th2} is now complete. 

\subsection{Proof of Theorem \protect\ref{nodetheorem}}\label
{subsecprthenode}

Under the sub-Gaussian assumption
we know that $\eta_j$ is also sub-Gaussian. So then $\| \eta_j^T \bx
_{-j} /n \|_{\infty}=
{\mathcal O}_{\PP} ( \sqrt{\log(p)/ n } )$. If $\| \bx\|_{\infty} =
{\mathcal O} (K) $, we can use the work in \cite{Duembgen09}
to conclude that
\[
\bigl\| \eta_j^T \bx_{-j} \bigr\|_{\infty}/n = {
\mathcal O }_{\PP} \bigl(K \sqrt {\log (p)/ n } \bigr).\vadjust{\goodbreak} %
\]
However, this result does not hold uniformly in $j$.
Otherwise, in the strongly bounded case, we have
\[
\| \eta_j \|_{\infty} \le\| X_j
\|_{\infty} + \| \bx_{-j} \gamma _j \|
_{\infty} ={ \mathcal O} (K). %
\]
So then $\| \eta_j^T \bx_{-j} /n\|_{\infty} = { \mathcal O }_{\PP} (K
\sqrt{\log(p)/ n } ) + {\mathcal O}_{\PP} ( K^2 \log(p)/ n ) $,
which is uniform in $j$.

Then by standard arguments (see, e.g., \cite{brt09}, and see \cite{pbvdg11}
which complements the concentration results in \cite{lederer11} for the
case of errors with only second moments) for $\lambda_j\asymp K_0
\sqrt
{\log(p)/ n } $ [recall that
$K_0 =1$ in the sub-Gaussian case and $K_0 = K$ in the (strongly)
bounded case]
\[
\bigl\| \bx_{-j} ( \hat\gamma_j - \gamma_{j} )
\bigr\|_n^2 = {\mathcal O}_{\PP} \bigl(s_j
\lambda_j^2 \bigr),\qquad \| \hat\gamma_j -
\gamma_{j} \|_1 = {\mathcal O}_{\PP} (
s_j \lambda_j ). %
\]
The condition $K^2 s_j \sqrt{ \log(p)/ n }$ is used in the (strongly)
bounded case to be able to conclude
that the empirical compatibility condition holds (see \cite{pbvdg11},
Section~6.12).
In the sub-Gaussian case, we use that $\sqrt{s_j \log(p)/n } = o(1)$
and an extension of Theorem 1 in \cite{raskutti2010restricted} from
the Gaussian case
to the sub-Gaussian case. This gives again that the empirical
compatibility condition holds.

We further find that
\[
\| \hat\gamma_j - \gamma_j \|_2 = {
\mathcal O}_{\PP} \bigl( K_0 \sqrt{s_j \log(p)/
n } \bigr). %
\]
To show this, we first introduce the notation $v^T \Sigma v:= \| \bx v
\|^2$. Then in the (strongly) bounded case
\[
\bigl| \| \bx v \|_n^2 - \| \bx v \|^2 \bigr| \le\|
\hat\Sigma- \Sigma\|_{\infty} \| v \|_1^2 = {
\mathcal O}_{\PP} \bigl( K^2 \sqrt{\log(p)/ n } \bigr) \| v
\|_1^2.
\]
Since $\| \hat\gamma_j - \gamma_j \|_1 = {\mathcal O}_{\PP} ( K_0s_j
\sqrt{\log(p)/ n } )$
and the smallest eigenvalue $\Lambda_\mathrm{min}^2 $ of $\Sigma$ stays away
from zero, this gives
\begin{eqnarray*}
{\mathcal O}_{\PP} \bigl( K_0^2
s_j \log(p)/ n\bigr)& =& \bigl\| \bx_{-j} ( \hat
\gamma_j - \gamma_j ) \bigr\|_n^2
\\
&\ge& \Lambda_\mathrm{min}^2 \| \hat\gamma_j -
\gamma_j \|_2^2 - {\mathcal O}_{\PP}
\bigl( K_0^4 s_j ^2 \bigl(
\log(p)/ n \bigr)^{3/2} \bigr)
\\
&\ge& \Lambda_\mathrm{min}^2 \| \hat\gamma_j -
\gamma_j \|_2^2 - o_{\PP}
\bigl(K_0^2 \log(p)/ n \bigr),
\end{eqnarray*}
where we again used that $K_0^2 s_j \sqrt{\log(p)/ n} = o(1)$.
In the sub-Gaussian case, the result for the $\| \cdot\|
_2$-estimation error follows by similar arguments invoking again
a sub-Gaussian extension of Theorem 1 in \cite{raskutti2010restricted}.

We moreover have
\begin{eqnarray*}
\bigl| \hat\tau_j^2 - \tau_j^2\bigr| &=
&\underbrace{\bigl| \eta_j^T \eta_j / n - \tau
_j^2\bigr| }_{I} + \underbrace{\bigl|
\eta_j^T \bx_{-j} (\hat\gamma_j -
\gamma_j )/n \bigr|}_{\mathit{II}}
\\
&&{}+ \underbrace{\bigl| \eta_j^T \bx_{-j}
\gamma_j /n\bigr| }_{\mathit{III}} + \underbrace {\bigl|(
\gamma_j)^T \bx_{-j}^T
\bx_{-j} (\hat\gamma_j - \gamma_j ) / n
\bigr|}_{\mathit{IV}}.
\end{eqnarray*}
Now, since we assume fourth moments of the errors,
\[
I = {\mathcal O}_{\PP} \bigl( K_0^2
n^{-1/2} \bigr). %
\]
Moreover,
\[
\mathit{II} = {\mathcal O }_{\PP} \bigl( K_0 \sqrt{\log(p)/ n}
\bigr) \| \hat\gamma_j - \gamma_j \|_1 = {
\mathcal O }_{\PP} \bigl( K_0^2 s_j
\log(p)/ n\bigr). %
\]
As for $\mathit{III}$, we have
\[
\mathit{III}= {\mathcal O}_{\PP} \bigl( K_0 \sqrt{\log(p)/ n} \bigr)
\| \gamma_j \|_1 = {\mathcal O}_{\PP} \bigl(
K_0 \sqrt{s_j\log(p)/ n} \bigr) %
\]
since $\| \gamma_j \|_1 \le\sqrt{s_j} \| \gamma_j \|_2 = {\mathcal O}
( \sqrt{s_j} ) $.
Finally, by the KKT conditions,
\[
\bigl\| \bx_{-j}^T \bx_{-j} ( \hat
\gamma_j - \gamma_j ) \bigr\|_{\infty} /n = {\mathcal
O}_{\PP} \bigl( K_0 \sqrt{\log(p)/ n } \bigr), %
\]
and hence
\[
\mathit{IV} = {\mathcal O}_{\PP} \bigl( K_0 \sqrt{\log(p)/ n }
\bigr) \| \gamma_j \|_1 = {\mathcal O}_{\PP}
\bigl( K_0 \sqrt{s_j \log(p)/ n } \bigr). %
\]
So now we have shown that
\[
\bigl| \hat\tau_j^2 - \tau_j^2\bigr| = {
\mathcal O}_{\PP} \bigl( K_0 \sqrt{s_j \log (p)/
n } \bigr). %
\]
Since $1/\tau_j^2 = {\mathcal O} (1)$, this implies that also
\[
1 / \hat\tau_j^2 - 1 / \tau_j^2
= {\mathcal O}_{\PP} \bigl( K_0 \sqrt{s_j
\log(p)/ n } \bigr). %
\]

We conclude that
\begin{eqnarray*}
\|\hat\Theta_j - \Theta_j \|_1 &=&\bigl\| \hat
C_j / \hat\tau_j^2 - C_j/
\tau_j^2 \bigr\|_1 \\
&\le& \underbrace{\| \hat
\gamma_j - \gamma_j \|_1 / \hat
\tau_j^2 }_{i} %
+ \underbrace{\| \gamma_j \|_1 \bigl( 1/ \hat
\tau_j^2 - 1 / \tau_j^2
\bigr)}_{\mathit{ii}},
\end{eqnarray*}
where
\[
i = {\mathcal O}_{\PP} \bigl( K_0 s_j \sqrt{
\log(p)/ n } \bigr) %
\]
since $\hat\tau_j^2$ is a consistent estimator of $\tau_j^2$ and $1/
\tau_j^2 = {\mathcal O} (1) $,
and also
\[
\mathit{ii}= {\mathcal O}_{\PP} \bigl( K_0 s_j \sqrt{
\log(p)/ n } \bigr), %
\]
since $\| \gamma_j \|_1 = {\mathcal O} (\sqrt{s_j}) $.

Recall that
\[
\| \hat\gamma_j - \gamma_j \|_2 = {
\mathcal O}_{\PP} \bigl( K_0 \sqrt{s_j \log(p)/
n } \bigr). %
\]
But then
\begin{eqnarray*}
\| \hat\Theta_j - \Theta_j \|_2 &\le&\| \hat
\gamma_j - \gamma_j \| _2 / \hat
\tau_j^2 + \| \gamma_j \|_2
\bigl( 1/ \hat\tau_j^2 - 1 / \tau_j^2
\bigr) %
\\
&=& {\mathcal O}_{\PP} \bigl(K_0 \sqrt{s_j
\log(p)/ n }\bigr). %
\end{eqnarray*}

For the last part, we write
\begin{eqnarray*}
&&\hat\Theta_j\Sigma\hat\Theta_j^T -
\Theta_{j,j} %
\\
&&\qquad= ( \hat\Theta_j - \Theta_j) \Sigma( \hat
\Theta_j - \Theta_j)^T +
\Theta_j \Sigma( \hat\Theta_j - \Theta_j
)^T + \Theta_j\Sigma\Theta_j^T
- \Theta_{j,j} %
\\
&&\qquad= ( \hat\Theta_j - \Theta_j) \Sigma( \hat
\Theta_j - \Theta_j)^T + 2\bigl(1/ \hat
\tau_j^2 - 1/ \tau_j^2 \bigr),
\end{eqnarray*}
since $\Theta_j \Sigma= e_j^T$, $\Theta_j\Sigma\Theta_j^T = \Theta_{j,j}$,
$\hat\Theta_{j,j} = 1/ \hat\tau_j^2$, and $\Theta_{j,j} = 1/ \tau
_j^2 $.
But
\[
( \hat\Theta_j - \Theta_j) \Sigma( \hat
\Theta_j - \Theta_j)^T \le\| \Sigma
\|_{\infty} \| \hat\Theta_j - \Theta_j
\|_1. %
\]
We may also use
\[
( \hat\Theta_j - \Theta_j) \Sigma( \hat
\Theta_j - \Theta_j)^T \le
\Lambda_\mathrm{max}^2 \| \hat\Theta_j -
\Theta_j \|_2^2.
\]

The last statement of the theorem follows as in Theorem \ref{th1},
as\break
$\sqrt n ( \hat b_{\mathrm{Lasso}, j} - \beta_j^0 ) = W_j + \Delta_j $, with
$\Delta_j \le\sqrt{n} \lambda_j / \hat\tau_j^2 \| \hat\beta-
\beta
^0 \|_1 $,
and $\lambda_j / \hat\tau_j^2 \asymp\lambda_j \asymp\sqrt{\log(p)
/ n }$, the latter being
uniformly in $j$ in the sub-Gaussian or strongly bounded case.

\subsection{Proof of Theorem \protect\ref{th3}}\label{th3section}

Note that
\[
\dot\rho( y, x_i \hat\beta) = \dot\rho\bigl( y, x_i
\beta^0\bigr) + \ddot\rho( y, \tilde a_i)
x_i\bigl( \hat\beta- \beta^0\bigr), %
\]
where $\tilde a_i$ is a point intermediating $x_i \hat\beta$ and $x_i
\beta^0$, so that
$| \tilde a_i - x_i \hat\beta| \le\break | x_i(\hat\beta- \beta^0) | $.

We find by the Lipschitz condition on $\ddot\rho$ [condition (C1)]
\begin{eqnarray*}
&&\bigl| \ddot\rho( y, \tilde a_i) x_i \bigl( \hat\beta-
\beta^0\bigr) - \ddot \rho ( y, x_i \hat\beta)
x_i \bigl( \hat\beta- \beta^0\bigr) \bigr| %
\\
&&\qquad\le | \tilde a_i - x_i \hat\beta| \bigl| x_i
\bigl(\hat\beta- \beta^0\bigr) \bigr| \le\bigl| x_i \bigl( \hat
\beta- \beta^0\bigr) \bigr|^2. %
\end{eqnarray*}
Thus, using that by condition (C5) $| x_i \hat\Theta_j^T | =
{\mathcal
O}_{\PP} (K) $ uniformly in $j$,
\[
\hat\Theta_j P_n \dot\rho_{\hat\beta} = \hat
\Theta_j P_n \dot \rho _{\beta^0} + \hat
\Theta_j P_n \ddot\rho_{\hat\beta} \bigl( \hat\beta-
\beta^0\bigr) + \mathrm{Rem}_1, %
\]
where
\begin{eqnarray*}
\mathrm{Rem}_1 &=& {\mathcal O}_{\PP} (K) \sum
_{i=1}^n \bigl| x_i \bigl( \hat\beta-
\beta^0\bigr) \bigr|^2/n = {\mathcal O}(K)\bigl \| \bx\bigl( \hat
\beta- \beta^0\bigr) \bigr\|_n^2 %
\\
&=& {\mathcal O}_{\PP} \bigl(Ks_0 \lambda^2
\bigr) = o_{\PP} (1), %
\end{eqnarray*}
where we used condition (C2) and in the last step condition (C8).

We know that by condition (C4)
\[
\bigl\| \hat\Theta_j P_n \ddot\rho_{\hat\beta} -
e_j^T \bigr\|_{\infty} = {\mathcal O} ( \lambda_*).
\]
It follows that
\begin{eqnarray*}
b_j- \beta_j^0 &=& \hat
\beta_j - \beta_j^0 - \hat
\Theta_j P_n \dot \rho_{\hat\beta}
\\
&=& \hat\beta_j - \beta_j^0 - \hat
\Theta_j P_n \dot\rho_{\beta
^0} - \hat
\Theta_j P_n \ddot\rho_{\hat\beta} \bigl( \hat\beta-
\beta^0\bigr) - \mathrm{Rem}_1
\\
&=& - \hat\Theta_j P_n \dot\rho_{\beta^0} - \bigl(
\hat\Theta_j P_n \ddot \rho_{\hat\beta} -
e_j^T\bigr) \bigl( \hat\beta- \beta^0\bigr)
- \mathrm{Rem}_1 \\
&=& -\hat\Theta_j P_n \dot
\rho_{\beta^0} - \mathrm{Rem }_2,
\end{eqnarray*}
where
\[
|\mathrm{Rem}_2 | \le|\mathrm{Rem}_1  | + {\mathcal O} (
\lambda_*)\bigl \| \hat \beta- \beta^0\bigr \|_1 =o_{\PP}
\bigl(n^{-1/2} \bigr) + {\mathcal O}_{\PP} (s_0
\lambda\lambda_*) = o_{\PP} \bigl(n^{-1/2} \bigr) %
\]
since by condition (C2) $\| \hat\beta- \beta^0 \|_1 = {\mathcal
O}_{\PP} ( \lambda s_0 ) $, and by
the second part of
condition~(C8) also $\lambda_* \lambda s_0 = o (n^{-1/2} ) $.

We now have to show that our estimator of the variance is consistent.
We find
\begin{eqnarray*}
& &\bigl|\bigl(\hat\Theta P \dot\rho_{ \beta^0 } \dot\rho_{\beta^0}^T
\hat{\Theta}^T \bigr)_{j,j} - \bigl( \hat\Theta
P_n \dot\rho_{\hat\beta} \dot \rho_{\hat\beta}^T
\hat\Theta^T \bigr)_{j,j} \bigr|
\\
&&\qquad\le \underbrace{\bigl| \bigl(\hat\Theta(P_n - P) \dot
\rho_{ \beta^0 } \dot\rho_{\beta^0}^T \hat
\Theta^T \bigr)_{j,j}\bigr |}_{I} + \underbrace{ \bigl|
\bigl(\hat\Theta P \dot\rho_{ \beta^0 } \dot\rho_{\beta^0}^T
\hat \Theta ^T \bigr)_{j,j} - \bigl( \hat\Theta P \dot
\rho_{\hat\beta} \dot\rho_{\hat\beta}^T \hat
\Theta^T \bigr)_{j,j} \bigr| }_{\mathit{II}}.
\end{eqnarray*}
But, writing $\varepsilon_{k,l}:= (P_n - P) \dot\rho_{k, \beta^0}
\dot\rho_{l, \beta^0} $, we see
that
\begin{eqnarray*}
I &=& \bigl|\bigl( \hat\Theta(P_n - P) \dot\rho_{ \beta^0 } \dot
\rho_{\beta
^0}^T \hat\Theta^T \bigr)_{j,j}\bigr|=
\biggl|\sum_{k,l} \hat\Theta_{j,k} \hat
\Theta_{j,l} \varepsilon_{k,l} \biggr| \le\| \hat
\Theta_j\|_1^2 \| \varepsilon
\|_{\infty}\\
& =& {\mathcal O}_{\PP} \bigl(s_* K^2\lambda
\bigr),
\end{eqnarray*}
where we used conditions (C5) and (C6).

Next, we will handle $\mathit{II}$. We have
\begin{eqnarray*}
\dot\rho_{\hat\beta} (y,x) \dot\rho_{\hat\beta}^T (y,x)- \dot
\rho _{\beta^0} (y,x) \dot\rho_{\beta^0}^T (y,x) &=& \bigl[
\dot\rho^2( y - x \hat\beta) - \dot\rho^2 \bigl(y - x
\beta^0 \bigr) \bigr] x^T x %
\\
&:=& v (y,x) x^T x, %
\end{eqnarray*}
with
\[
\bigl| v(y,x)\bigr |:=\bigl | \dot\rho^2( y - x\hat\beta) - \dot\rho^2
\bigl(y- x\beta ^0 \bigr)\bigr| = {\mathcal O}_{\PP} (1)\bigl | x
\bigl( \hat\beta- \beta^0\bigr) \bigr|,
\]
where we use that $\dot\rho_{\beta^0} $ is bounded and $\ddot\rho$ is
locally bounded [condition (C1)].
It follows from condition (C2) that
\[
P |v | \le\sqrt{P |v|^2} =\bigl\| \bx\bigl( \hat\beta-
\beta^0 \bigr) \bigr\| = {\mathcal O}_{\PP} ( \lambda
\sqrt{s_0}). %
\]
Moreover, by condition (C5),
\[
\bigl\| \hat\Theta_j x^T \bigr\|_{\infty} = {\mathcal
O}_{\PP} (K ) %
\]
so that
\[
\bigl| \bigl( \hat\Theta v(x,y) x^T x \hat\Theta^T
\bigr)_{j,j}\bigr | \le{\mathcal O} \bigl(K^2 \bigr) \bigl|v(y,x)\bigr |.\vadjust{\goodbreak}
\]
Thus,
\[
\bigl| \bigl(\hat\Theta P \dot\rho_{ \beta^0 } \dot\rho_{\beta^0}^T
\hat \Theta ^T \bigr)_{j,j} - \bigl( \hat\Theta P \dot
\rho_{\hat\beta} \dot\rho_{\hat\beta}^T \hat
\Theta^T\bigr)_{j,j}\bigr |= {\mathcal O}_{\PP} \bigl(
K^2 \sqrt{s_0} \lambda\bigr). %
\]
It follows that
\[
I + \mathit{II} = {\mathcal O}_{\PP} \bigl( K^2 s_* \lambda\bigr) +
{\mathcal O}_{\PP} \bigl( K^2 \sqrt{s_0} \lambda
\bigr) = o_{\PP} (1) %
\]
by the last part of condition (C8). 

\subsection{Proof of Theorem \protect\ref{GLMtheorem}}\label
{GLMtheoremsection}
This follows from Theorem \ref{th3}.
The assumptions~(C2), (C4)--(C8)
follow from the conditions of Corollary \ref{node2corollary} with
$\Sigma_{\beta}:= P \ddot\rho_{\beta}$ and
$w_{\beta}^2 (y,x):= \ddot\rho(y, x \beta) $,
where we take $\hat\Theta= \hat\Theta_{\mathrm{Lasso}}$
and $s_* = s_j$ and $\lambda_* = \lambda_j $.
Condition (C2) holds because the compatibility condition is met as
$\Sigma_{\beta^0}$ is nonsingular and
\[
\| \hat\Sigma- \Sigma_{\beta^0} \|_{\infty} = { \mathcal
O}_{\PP} ( \lambda_*). %
\]
The condition that
$\dot\rho(y, x \beta^0) $ is bounded ensures that $\rho(y,a)$ is
locally Lipschitz, so that we can control the empirical process
$(P_n- P) ( \rho_{\hat\beta} - \rho_{\beta^0} ) $ as in \cite
{geer07} (see also~\cite{pbvdg11} or \cite{Mueller11}). [In the case of a GLM with
canonical loss
(e.g., least squares loss) we can relax the
condition of a locally bounded derivative because the empirical process
is then linear.]
Condition (C3) is assumed to hold with $\|\bx\|_{\infty} = {\mathcal O}
(1) $, and condition (C4) holds with $\lambda_* \asymp\sqrt{\log p/
n} $.
This is because
in the nodewise regression construction,
the $1/ \hat\tau_j^2$ are consistent estimators of $(\Sigma_{\beta
^0}^{-1} )_{jj}$
(see Theorem \ref{node2theorem}).
Condition (C5) holds as well. Indeed, $ \| \Theta_{\beta^0,j }\|_1 =
{\mathcal O} (\sqrt{s_j})$, and
$\| \hat\Theta_{\hat{\beta},j}- \Theta_{\beta^0, j} \|_1 =
{\mathcal
O}_{\PP}(
\lambda_j s_j) = {\mathcal O}_{\PP} ( \sqrt{s_j}) $.
Condition (C6) holds, too, since we assume that $\| \dot\rho_{\beta
^0} \|_{\infty} = {\mathcal O} (1)$ as well as
$\|\bx\|_{\infty} = {\mathcal O} (1)$.
As for condition (C7), this follows from Lemma \ref{normallemma}, since
$| \Theta_{\beta^0,j}\dot\rho_{\beta^0} (y,x) | =| \Theta_{\beta^0,j}
x^T \dot\rho(y, x \beta^0) | = {\mathcal O} (1)$, which implies
for $A:= P \dot\rho_{\beta^0} \dot\rho_{\beta^0}^T $ that
$ \| A \Theta_{\beta^0,j}^T \|_{\infty} ={\mathcal O} (1) $. 

\begin{supplement}[id=suppA]
\stitle{Supplement to ``On asymptotically optimal confidence regions and
tests for high-dimensional models''}
\slink[doi]{10.1214/14-AOS1221SUPP} 
\sdatatype{.pdf}
\sfilename{aos1221\_supp.pdf}
\sdescription{The supplemental article contains additional empirical
results, as well as the proofs of Theorems
\ref{thspj} and \ref{node2theorem}, Lemmas \ref{concretelem} and
\ref{normallemma}.}
\end{supplement}

%


\printaddresses
\end{document}